\documentclass{article}

\usepackage{dsfont}                     

\usepackage{enumerate}                  

\usepackage{verbatim}

\usepackage{url}
\usepackage[colorlinks=true,urlcolor=black]{hyperref}                   

\usepackage{array, multirow, tabularx}  

\usepackage{soul}                       
\usepackage[normalem]{ulem}             

\usepackage{amsmath,amssymb,mathrsfs}   

\usepackage{graphicx}                   
\usepackage{float, caption}             
\usepackage{subfig}                     

\usepackage{calc,ifthen,xspace}         

\usepackage{hyperref}

\newcommand{\dem}[1]{\noindent \textbf{Proof:} #1 \hfill $\Box$}

\newcommand{\R}{\ensuremath{\mathbb{R}}}       

\newtheorem{Th}{Theorem}[section]
\newtheorem{Prop}[Th]{Proposition}
\newtheorem{Lemma}[Th]{Lemma}
\newtheorem{Cor}[Th]{Corollary}
\newtheorem{Conjecture}[Th]{Conjecture}

\newtheorem{Exa}[Th]{Example}
\newtheorem{Rem}[Th]{Remark}

\begin{document}
\thispagestyle{plain}
\begin{center}
  \Large
    On the sharpness of the weighted Bernstein-Walsh inequality, with applications to the superlinear convergence of conjugate gradients
  \\[20pt]
  \normalsize by \\[15pt]
  \large Bernhard Beckermann\footnote{Laboratoire Painlev\'e UMR 8524, UFR Math\'ematiques, Univ.\ Lille, F-59655 Villeneuve d'Ascq CEDEX, France. E-mail:
{\tt bbecker@math.univ-lille1.fr}. Supported  in  part  by  the  Labex  CEMPI
   (ANR-11-LABX-0007-01).} and Thomas Helart\footnote{Laboratoire Painlev\'e UMR 8524, UFR Math\'ematiques, Univ.\ Lille, F-59655 Villeneuve d'Ascq CEDEX, France. E-mail:
{\tt thomas.helart@gmail.com}.}
\end{center}




\def\supp{\rm supp}


\begin{abstract}
  In this paper we show that the weighted Bernstein-Walsh inequality in logarithmic potential theory is sharp up to some new universal constant, provided that the external field is given by a logarithmic potential. 
  Our main tool for such results is a new technique of discretization of logarithmic potentials, where we take the same starting point as in earlier work of Totik and of Levin \& Lubinsky, but add an important new ingredient, namely some new mean value property for the cumulative distribution function of the underlying measure.

  As an application, 
  we revisit the work of Beckermann \& Kuijlaars on the superlinear convergence of conjugate gradients. These authors have determined the asymptotic convergence factor for sequences of systems of linear equations with an asymptotic eigenvalue distribution. 
  There was some numerical evidence to let conjecture that the integral mean of Green functions occurring in their work
  should also allow to give inequalities for the rate of convergence if one makes a suitable link between measures and the eigenvalues of a single matrix of coefficients. 
  We prove this conjecture, at least for a class of measures which is of particular interest for applications.
\end{abstract}




\section{Introduction}

Conjugate gradients (CG) is a popular method for solving large systems of equations, with the matrix of coefficients being symmetric and positive definite \cite{Bau:Trefethen}.
However, its convergence (at least in exact arithmetic) is not yet fully understood, despite an important number of research contributions,
for instance \cite{Beck-CG,Beck-sharpness,Beck-CG:right-hand-side,Beck_Springer,
Driscoll-Trefethen:from-potential-theory, Ipsen, Kuijlaars-which-eigenvalues-2000, Kuijlaars:review, van-der-Sluis:Van-der-Vorst:86}.
It happens quite often that there is a regime of convergence called superlinear convergence, which depends very much on the eigenvalue distribution of the matrix, see \S\ref{sec_superlinear} for more details.
People have been aware of this phenomenon for more than 40 years, but only in \cite{Beck-CG} a general theory based on logarithmic potential theory was suggested to quantify the rate of convergence, see also \cite{Kuijlaars:review} for a more comprehensive summary. The drawback of this theory is that all results in \cite{Beck-CG} study only the so-called asymptotic convergence factor. In addition, this theory requires to consider sequences of systems of equations with a joint eigenvalue distribution, and thus gives not so much information about the actual rate of convergence for a single matrix.
Numerical evidence in \cite{Beck-CG,Beck-sharpness,Beck-CG:right-hand-side} seemed to indicate that behind the asymptotic results there should be some hidden inequality valid for a single matrix, see Conjecture~\ref{conjecture_superlinear} below. To our knowledge, the present paper is the first which deals with this conjecture, at least for a suitable subclass of eigenvalue distributions. 

This paper 
contains three main ingredients, all being connected with polynomial extremal problems and thus with logarithmic potential theory: we discuss in \S\ref{sec_BW} the sharpness of the so-called weighted Bernstein-Walsh inequality for the particular case where the external field is the logarithmic potential of some measure. Here our main result in Theorem~\ref{th_main} indicates the existence of some new universal constant. Secondly, we give and discuss in \S\ref{sec_superlinear} some new upper bound for the rate of convergence of conjugate gradients, and show in our Theorem~\ref{Th_super} the above conjecture for a particular class of eigenvalue distributions, which is illustrated by some (academic) numerical examples.

Our main technical result stated and proved in \S\ref{sec_discretization} is Theorem~\ref{theorem-potential} on a new fine discretization of logarithmic potentials for a suitable class of measures, where in contrast to preceding work of Totik, Lubinsky and others we get (large but) explicit constants. Here an essential tool is a new mean value property stated in Theorem~\ref{Theorem:mean:property}.

\subsection{The weighted Bernstein-Walsh inequality}\label{sec_BW}

One of the appealing aspects of CG convergence is that there is a close link with polynomial extremal problems and extremal problems in logarithmic potential theory, which we discuss now.

Given a finite union of compact intervals $\Sigma\subset \R$, we denote by $\mathcal{M}_1(\Sigma)$ the set of Borel measures $\mu$ with support
$\supp(\mu)$ in $\Sigma$ and of total mass $1$, and consider the logarithmic potential and energy
\[
   U^{\mu}(x) = \int \log \frac{1}{|x-t|} \mathrm{d}\mu(t) , \quad
   I(\mu) = \int \int \log \frac{1}{|x-t|} \mathrm{d}\mu(t) \mathrm{d}\mu(x).
\]
Given a weight $w$ defined on $\R$ and continuous on $\Sigma$ together with an external field  $Q(x)=-\log (w(x))$,
it is known \cite[Theorem I.1.3 and Theorem I.4.8]{Saff-Totik} that there is a unique minimizer
$\mu\in \mathcal M_1(\Sigma)$ of the extremal problem
\begin{equation} \label{eq.extremal1}
    \inf \{ I(\mu) + 2 \int Q \mathrm{d}\mu : \mu \in \mathcal M_1(\Sigma) \}
\end{equation}
which is uniquely characterized by the existence of a constant $F\in \mathbb R$ such that
\begin{equation} \label{eq.extremal2}
    f(x):=F - U^\mu(x) - Q(x) \left\{\begin{array}{ll} = 0  & \mbox{for $x\in \supp(\mu)$},
    \\
    \leq 0  & \mbox{for $x\in \Sigma$.}
\end{array}\right.
\end{equation}

Logarithmic potential theory with external fields has been applied with success for getting asymptotics for various polynomial extremal problems \cite{Saff-Totik},
maybe one of the most prominent results being the weighted Bernstein-Walsh inequality \cite[Theorem~III.2.1]{Saff-Totik}
\begin{equation} \label{eq.extremal3}
      \forall x_0\in \mathbb R \quad \forall P \in \Pi_k:
      \quad
      \frac{| w(x_0)^k P(x_0) |}{\| w^k P \|_{\supp(\mu)}}  \leq e^{kf(x_0)} ,
\end{equation}
and its sharpness, see, e.g., \cite[Corollary~III.1.10]{Saff-Totik},
\begin{equation} \label{eq.extremal4}
      \exists P_k \in \Pi_k \quad \forall x_0\in \mathbb R \setminus \supp(\mu) :  \quad
      \lim_{k\to \infty} \left( \frac{| w(x_0)^k P_k(x_0) |}{\| w^k P_k \|_{\supp(\mu)}} \right)^{1/k} = e^{f(x_0)} ,
\end{equation}
where $\Pi_k$ denotes the set of polynomials of degree at most $k$, and $\| f \|_{\Sigma} = \max_{x\in \Sigma}|f(x)|$. One aim of this paper is to improve
\eqref{eq.extremal4} for a particular class of external fields, see Theorem~\ref{th_main} below,
namely to show that \eqref{eq.extremal3} is sharp up to some constant. Before giving some more details,
let us first have a look at other classes of external fields where such constants are explicitly known.
In what follows we will write $g_S(\cdot,\zeta)$ to denote the Green function in $\overline{\mathbb C} \setminus S$ for a compact set $S\subset \mathbb R$ with pole at $\zeta\in \overline{\mathbb C} \setminus S$. We will be mainly interested in the special case of an interval $S$ where the Green function vanishes on $S$ and is strictly positive outside $S$, and where explicit formulas are available.

\begin{Exa}\label{example1}
   Consider $\Sigma=[a,b]$ and $Q=0$, then an explicit formula is known for the minimizer in \eqref{eq.extremal1}
   $$
        \supp(\mu)=\Sigma=[a,b], \quad
        \frac{d\mu}{dx}(x) = \frac{1}{\pi \sqrt{(x-a)(b-x)}},
   $$
   also called Robin equilibrium measure of the interval $[a,b]$ and denoted by $\omega_{[a,b]}$. It is also known from, e.g., \cite[Eqn.\ (I.4.8)]{Saff-Totik} that
   $f(z)=g_{[a,b]}(z,\infty)$, and thus
   \eqref{eq.extremal3} becomes the classical Bernstein-Walsh inequality.
	Taking $P_k(x)=T_k(\frac{2x-a-b}{b-a})$ with $T_k$ the Chebyshev polynomial of the first kind, one may also show that \eqref{eq.extremal3} is sharp up to a factor $1/2$.
\end{Exa}

\begin{Exa}
   Consider $\Sigma=[a,b]$, and $w(x)^k=1/\sqrt{q(x)}$ with $q$ being a polynomial of degree $\ell\leq 2k$, strictly positive on $[a,b]$,
	compare with \cite[chap 4.4]{Meinardus}. Thus $Q=-U^\rho$ with $\rho $ an atomic measure of mass $\ell/(2k)\leq 1$.
	Here the extremal measure in \eqref{eq.extremal1}, \eqref{eq.extremal2} is given in \cite[Example~II.4.8]{Saff-Totik}
	in terms of balayage onto $\supp(\mu)=\Sigma$, and it follows from \cite[Eqn.\ (4.32)]{Saff-Totik} that
   \[
        f(x) = (1-\frac{\ell}{2k}) g(x,\infty) + \int g(x,y) \mathrm{d}\rho(y) .
   \]
   Moreover, with help of the factorization
   \[
         \widetilde q(y) \widetilde q(\frac{1}{y})
         =  q(x) , \quad \frac{2x-a-b}{b-a}= \frac{1}{2} (y+\frac{1}{y}) \in \mathbb R ,
   \]
   $|y| \geq 1$,
   the polynomial $\widetilde q$ of degree $\ell$ having all its roots outside the unit circle, it is known that $P_k$ defined by
   \[
        w(x)^k P_k(x) = \frac{1}{2} ( e^{kf(x)} + e^{-kf(x)}) , \quad
        e^{2kf(x)} = \frac{y^{2k} \widetilde q(\frac{1}{y})}{\widetilde q(y)},
   \]
   is a polynomial of degree $k$, showing that again \eqref{eq.extremal3} is sharp up to a factor $1/2$. 
\end{Exa}

We are interested in the case where the external field is a positive potential $U^{\rho/k}$ (not necessarily of an atomic measure), for instance if $w^k$ is a (power of a) polynomial. This includes
the particular case $w(x)=|x|^\theta$ on $\Sigma=[0,1]$ for $\theta>0$, starting point of an important research area about incomplete polynomials \cite[\S VI.1.1]{Saff-Totik}. For external fields being a positive potential, we recall below how to solve the extremal problem, including the well-known pushing effect that the support of the equilibrium measure may be a proper subset of $\Sigma$. We then state our main result on the sharpness of the weighted Bernstein-Walsh inequality.

\begin{Th}\label{th_main}
   Let $k \geq 1$ be some integer, and $Q=U^{\rho/k}$ on $\Sigma=[\alpha,b]$, with the Borel measure $\rho$ being compactly supported on $(-\infty,\alpha]$. Consider on $\Sigma$ the strictly decreasing function
   \begin{equation}
      \label{condition:mu:positive}
         \eta(a):= \int \sqrt{\frac{b-y}{a-y}} \mathrm{d} \rho(y),
   \end{equation}
   and set $a=\alpha$ if $k + || \rho || \geq \eta(\alpha)$, and else denote by $a$ the unique solution of $k + || \rho || = \eta(a)$ in $\Sigma$.
   Then the extremal measure in \eqref{eq.extremal1}, \eqref{eq.extremal2} is given by
   \begin{equation} \label{eq.extremal5}
        \supp(\mu)=[a,b], \quad k f(x) = (k+\| \rho\|) g_{[a,b]}(x,\infty) -
        \int g_{[a,b]}(x,y) \mathrm{d}\rho(y) .
   \end{equation}
   Moreover, the weighted Bernstein-Walsh inequality \eqref{eq.extremal3} is sharp up to some constant, that is,
	there exists a universal real constant $C_{BW}>0$ such that, for all $k \geq 2$, 
	we may construct a polynomial $P_k$ of degree $k$ such that, for all $x_0\in \mathbb R \setminus [a,b]$,
   \begin{equation} \label{eq.extremal6}
      \frac{| w(x_0)^k P_k(x_0) |}{\| w^k P_k \|_{\supp(\mu)}}  \geq e^{-C_{BW}} \,   e^{kf(x_0)} .
   \end{equation}
\end{Th}

Our proof of Theorem~\ref{th_main} presented in \S\ref{subsec_proof_th_main} is based on a fine discretization of the logarithmic potential $U^{k\mu}$.
We will show in this paper that $C_{BW}\leq 15 383$, but this is by no means optimal.
The most remarkable fact for us seems to be that such a constant does not depend on the data $\rho,a,b$ nor on $k$.
In particular, we do not need any further assumptions on smoothness of $\rho$, which is probably required by other techniques like a Riemann-Hilbert approach
(which in any case would only allow to discuss asymptotics).

\subsection{Superlinear convergence for conjugate gradients}\label{sec_superlinear}

\begin{figure}[t]
   \centerline{\includegraphics[width=0.7\textwidth,height=0.49\textwidth]{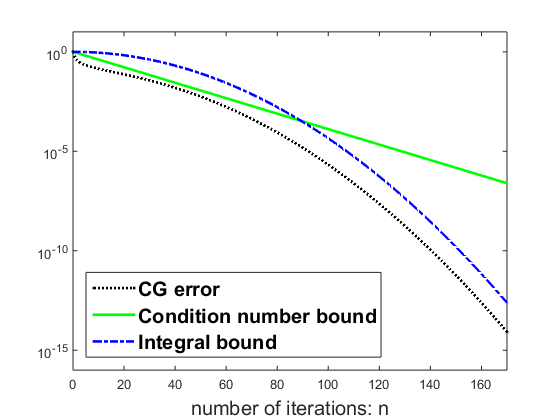}}
   \caption{Lower and upper bounds for $n  \mapsto E_n(\Lambda(A))$. Here
   $\lambda_j=(j/N)\, (2-j/N)$ for $j=1,...,N=1000$. As lower bound we draw the relative CG error in energy norm,
	 with $A=\mbox{diag}(\lambda_1,...,\lambda_N)$, $c=(1,...,1)^T$, and starting vector $x_0^{CG}=0$ (black dotted line).
	 The upper bounds come from \eqref{CG_bound} with $b/a=\lambda_N/\lambda_1$ for the condition number bound (green solid line)
	 and from Conjecture~\ref{conjecture_superlinear} with $C=0$ (blue dash-dotted line), see Example~\ref{example11} for further details.
	 }\label{fig1}
\end{figure}

Conjugate gradients is a popular method for solving large sparse linear systems $A x = c$ with symmetric positive definite $A$,
with spectrum $\Lambda(A)=\{ \lambda_j \}$, $0<\lambda_1<\lambda_2<...\leq b$.
Here one easily obtains the error estimate\footnote{In general, \eqref{CG_error} might be an important overestimation of the error, but there exist right-hand sides $c$ with equality.}
for the $n$th iterate $x_n^{CG}$
\begin{equation} \label{CG_error}
    \frac{\| x_n^{CG} - A^{-1} c \|_A}{\| x_0^{CG} - A^{-1} c  \|_A}
    \leq \min_{q\in \Pi_n} \max_{x\in S} \left|\frac{q(x)}{q(0)} \right| =: E_n(S)
\end{equation}
with the energy norm $\| y \|_A^2=y^*Ay$, where $S$ is any compact set containing the spectrum $\Lambda(A)$, for instance $S=\Lambda(A)$.
Thus, in contrast to the polynomial extremal problem considered in \S\ref{sec_BW}, we have a trivial weight and take norms on discrete sets.
One way of relating the two problems is to replace $\Lambda(A)$ by an interval $[a,b]\subset (0,\infty)$ containing all eigenvalues,
leading to the classical upper bound 
\begin{equation} \label{CG_bound}
      E_n(\Lambda(A)) \leq 2 \, \exp (-ng_{[a,b]}(0,\infty)) = 2 \, \Bigl( \frac{\sqrt{b/a}-1}{\sqrt{b/a}+1}\Bigr)^n,
\end{equation}
compare with Example~\ref{example1}. It is however known for a long time that there are eigenvalue distributions which lead to convergence which is faster than the one described in \eqref{CG_bound}, namely so-called superlinear convergence, see for instance Figure~\ref{fig1}.
A first attempt to quantify such a convergence behavior was suggested by Kuijlaars and Beckermann \cite{Beck-CG}, see also the 
review \cite{Kuijlaars:review}, or the review \cite{Beck_Springer} from the perspective of discrete orthogonal polynomials.
The key ingredient of this theory is to dispose of a measure $\sigma$ with continuous potential $U^\sigma$ and compact support describing the eigenvalue distribution.
In \cite{Beck-CG}, this is quantified by supposing that there is a sequence of systems $A_Nx_N=c_N$, with $\sigma$ being
the weak-star limit of normalized counting measures of the spectra of the symmetric and positive definite matrices $A_N$,
\begin{equation} \label{weak_limit}
      \lim_{N\to \infty} \frac{1}{N} \sum_{\lambda\in \Lambda(A_N)} \delta_\lambda = \sigma. 
\end{equation}
Under some additional weak assumptions for small eigenvalues,
the authors establish in \cite[Theorem~2.1]{Beck-CG} for the $n$th iterate $x_{n,N}^{CG}$ of conjugate gradients applied to the system $A_Nx_N=c_N$ the asymptotic upper bound
\begin{eqnarray}\nonumber
\limsup_{{n,N\to \infty}\atop {n/N \to t}} \left( \frac{\| x_{n,N}^{CG} - A_N^{-1} c_N \|_{A_N}}{\| x_{0,N}^{CG} - A_N^{-1} c_N  \|_{A_N}} \right)^{1/n}
    &\leq&
    \limsup_{{n,N\to \infty}\atop {n/N \to t}} \left( E_n(\Lambda(A_N))\right)^{1/n}
    \\&\leq&\label{eq.BB_ABK}
    \exp\left(-\frac{1}{t}\int_0^t g_{S(\tau)}(0,\infty) d\tau \right),
\end{eqnarray}
where $(S(t))_{0<t<\|\sigma\|}$ is a decreasing family of compact subsets of the convex hull of the spectra, obtained from some constrained extremal problem in logarithmic potential theory, which we explain now.

For measures $\sigma$ with compact support and continuous potential, and $0<t<\|\sigma\|$, according to \cite{Dragnev-Saff-1997, Rak96}
there exists a unique minimizer $\nu_{t,\sigma}$ of $I(\nu)$ under all candidates $\nu\in \mathcal M_1(\supp(\sigma))$ with $\nu \leq \sigma/t$.
This minimizer is uniquely characterized by the existence of a constant $C_{t,\sigma}\in \mathbb R$ such that
$$
     U^{\nu_{t,\sigma}}(x) = C_{t,\sigma} \mbox{~~for $x\in \supp(\sigma/t-\nu_{t,\sigma})$}, \quad
     U^{\nu_{t,\sigma}}(x) \leq C_{t,\sigma} \mbox{~~for $x\in \supp(\sigma)$}.
$$
Many Buyarov-Rakhmanov type properties are known about the measures $\nu_{t,\sigma}$ for fixed $\sigma$ and varying $t$, we just recall here from \cite[Proof of Theorem~2.1]{Beck-CG} the fact
that the measures $t\nu_{t,\sigma}$ are increasing in $t$, and hence
$$
     S(t):=\supp(\sigma/t-\nu_{t,\sigma}) \quad \mbox{decreases in $t$.}
$$
As a consequence, the map $n \mapsto - N \int_0^{n/N} g_{S(\tau)}(0,\infty) d\tau$ is concave and describes superlinear convergence behavior.
The compact sets $S(t)$ may have a quite complicated shape, and the main finding of \cite{Kuijlaars-which-eigenvalues-2000} roughly says that
the $n$th Ritz values of $A_N$ approach well all eigenvalues in $\Lambda(A_N) \setminus S(n/N)$. There is a similar (rough) interpretation of \eqref{eq.BB_ABK}:
so-called "converged" eigenvalues which are already well approached by $n$th Ritz values should no longer contribute (in exact arithmetic) to the convergence of CG at later stages.

In many examples, numerical evidence did let to conjecture that the above upper bound \eqref{eq.BB_ABK} even holds (up to some modest constant) for a single matrix $A$,
without limits and without taking the $n$-th root, see for instance \cite[Eqn.(1.9) and Figures 1 and 4]{Beck-CG}, \cite[Eqn.\ (1.3)]{Beck-CG:right-hand-side}, or Figure~\ref{fig1}.
Of course, for a single matrix we cannot define $\sigma$ through \eqref{weak_limit}. This gives the following conjecture.

\begin{Conjecture}\label{conjecture_superlinear}
   There is a (modest) constant $C\in \mathbb R$ and a technique of associating a measure $\sigma$ with compact support and continuous potential to the spectrum of a positive definite matrix $A$ such that, for all $n$ sufficiently small,
   $$
      E_n(\Lambda(A)) \leq \exp\left( C - N \int_0^{n/N} g_{S(t)}(0,\infty) dt \right).
   $$
\end{Conjecture}

It may be that this conjecture is wrong for measures where $S(t)$ has a complicated shape.
In our proof of the conjecture, following \cite[Lemma 3.1(a)]{Beck-CG}, we will impose sufficient conditions on $\sigma$ such that $S(t)=[a(t),b]$ for all $t$.

\begin{Lemma}\label{lemma_CG}
  Suppose that $\sigma$ is supported on the interval $[a(0),b]$ with density with respect to Lebesgue measure denoted by $\sigma'$, and suppose\footnote{It follows that $\sigma$ has compact support and continuous potential.} that $x\mapsto \sqrt{(x-a(0))(b-x)} \sigma'(x)$ vanishes at $x=a$, and is strictly increasing in $[a(0),b]$. Then for all $t\in (0,\| \sigma\|)$ we have $S(t)=[a(t),b]$, with $a(t)$ being the unique solution of the equation
  $$
         t = \int_{a(0)}^{a(t)} \sqrt{\frac{b-x}{a(t)-x}} d\sigma(x) ,
  $$
  in particular $t\mapsto a(t)$ is strictly increasing.
\end{Lemma}

Roughly speaking, having $S(t)=[a(t),b]$ for sufficiently small $t$ means that there are so few eigenvalues around $0$ that they are the first eigenvalues which are well approached by Ritz values of low order.
One of the reasons to consider such sets $S(t)$ is that, in any case, the superlinear convergence rate is only pronounced if small eigenvalues are well approached by Ritz values,
and the rate depends not as much on other "converging" eigenvalues, which in first order could be neglected.
Another reason is that, if the system $Ax=c$ comes from discretizing an elliptic PDE,
we might have only asymptotic knowledge on small eigenvalues of $A$ through a so-called Weyl formula.
The final reason 
is that in the particular case $S(t)=[a(t),b]$ the analysis becomes simpler, and also the upper bound is more explicit, since, by \eqref{CG_bound},
\begin{eqnarray} \nonumber
    \exp\Bigl( - N \int_0^{n/N} g_{S(t)}(0,\infty) dt \Bigr)
    &=&
    \exp\Bigl( N \int_0^{n/N} \log(\frac{\sqrt{b/a(t)}-1}{\sqrt{b/a(t)}+1})  dt \Bigr)
    \\&\leq& \prod_{j=0}^{n-1} \frac{\sqrt{b/a(j/N)}-1}{\sqrt{b/a(j/N)}+1}
\end{eqnarray}
in terms of some "effective condition number" $b/a(j/N)$, compare with \cite[Eqn.\ (2.27)]{Beck-CG:right-hand-side}.

In order to proceed, we first extend our definition \eqref{CG_error} of $E_n(S)$ to compact sets $S$ which do not necessarily contain the spectrum of $A$: following \cite{Beck-sharpness},
for a fixed matrix $A$, a compact set $S$, and sufficiently large $n$, let
$$
     E_n(S):= \min_{q\in \Pi_n} \left\{ \frac{\| q \|_{S}}{|q(0)|}: \forall \lambda\in \Lambda(A)\setminus S, q(\lambda)=0 \right\},
$$
and then obviously $E_n(\Lambda(A)) \leq E_n(S)$.
This inequality has been used for example in \cite{Ipsen} or \cite{van-der-Sluis:Van-der-Vorst:86} in order to derive a CG convergence bound
taking into account few outliers represented by the set $\Lambda(A)\setminus S$,
where typically $S$ is the convex hull of the remaining eigenvalues.

In what follows we consider $S=[\lambda_{d+1},b]$, and thus we prescribe
as roots of $q$ the smallest $d$ eigenvalues $\lambda_{1},...,\lambda_d$.
Understanding the modulus of the product of the corresponding linear factors as a weight,
and setting $\rho=\delta_{\lambda_1}+...+\delta_{\lambda_d}$, $\alpha=\lambda_{d+1}$ and $n=k+d=k+\|\rho\|$, Theorem~\ref{th_main} gives the following upper bounds in terms of Green functions. The sharpness follows from the weighted Bernstein-Walsh inequality \eqref{eq.extremal3}.

\begin{Cor}\label{Cor_bound}
   For any integer $n>d+1\geq 1$, let $a=a_{d,n}$ be equal to $\lambda_{d+1}$ if $n \geq \sum_{j=1}^d \sqrt{\frac{b-\lambda_j}{\lambda_{d+1}-\lambda_j}}=\eta(\lambda_{d+1})$, and else let $a$ be the unique solution $>\lambda_{d+1}$ of the equation $n=\sum_{j=1}^d \sqrt{\frac{b-\lambda_j}{a-\lambda_j}}$. Then
   \begin{equation} \label{eq.End}
     E_{n}(([\lambda_{d+1},b]) \leq \exp\Bigl(C_{BW} - n g_{[a,b]}(0,\infty)+ \sum_{j=1}^d g_{[a,b]}(0,\lambda_j)\Bigr),
   \end{equation}
   being sharp up to the factor $\exp(C_{BW})$.
\end{Cor}

\begin{figure}[t]
   \centerline{\includegraphics[width=0.49\textwidth]{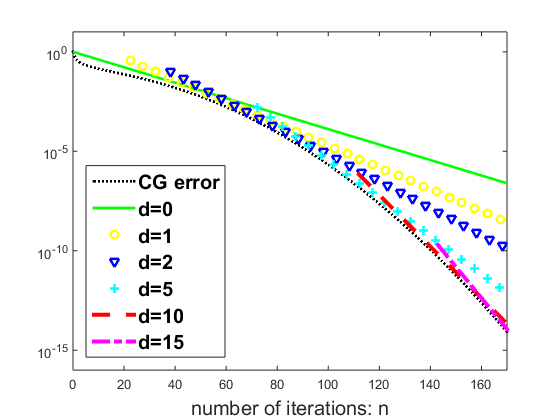}~
   \includegraphics[width=0.49\textwidth]{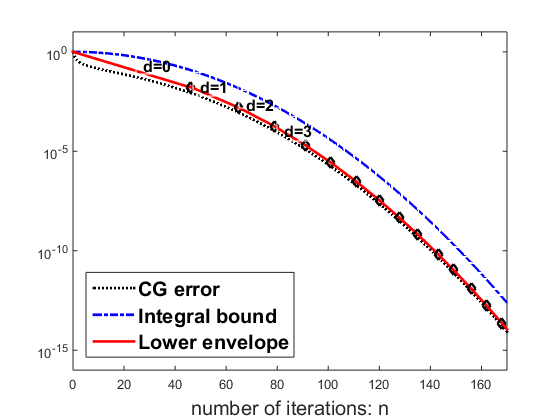}}
   \caption{Illustration of Corollary~\ref{Cor_bound}.
   As lower bound we draw on both plots the relative CG error in energy norm, with $\lambda_j,N,A,c,x_0^{CG}$ as in Figure~\ref{fig1} (black dotted line).
	 The straight lines on the left correspond to the bounds for $d\in \{ 0,1,2,5,10,15\}$ given in Corollary~\ref{Cor_bound},
	 each time for $n$ sufficiently large such that $a_{n,d}=\lambda_{d+1}$. Since it is difficult to see details,
	 we have drawn on the right only the CG error and the polygon obtained from the lower envelope of all straight lines in Corollary~\ref{Cor_bound},
	 where we indicate in the plot the correspondence between a segment and the choice of $d$.
	 To compare, we also have reproduced on the right from Figure~\ref{fig1} the integral bound from Conjecture~\ref{conjecture_superlinear}
	with $C=0$ (blue dash-dotted line),
	 verifying numerically that Corollary~\ref{Cor_bound} is the right tool to prove the conjecture.
	 }\label{fig2}
\end{figure}
Corollary~\ref{Cor_bound} gives us for each $d<n-1$ an upper bound for the function
$n \mapsto \log E_n(\Lambda(A))$,
each of them having the shape of a straight line for sufficiently large $n$, with the slope $-g_{[\lambda_{d+1},b]}(0,\infty)$
of these straight lines decreasing with $d$, but the abscissa in general increases. We thus hope that $\log E_n(\Lambda(A))$
is close to the value of the concave lower envelope of these straight lines, which is true for the particular example of Figure~\ref{fig2}.
In fact, finding an optimal $d=d_n<n$ with minimal $E_{n}(([\lambda_{d+1},b])$ for given $n$ seems to be a difficult task,
we will suggest an approximate solution in order to solve the above conjecture.

\begin{Th}\label{Th_super}
  Let $\sigma$ and $S(t)=[a(t),b]$ for $0<t<\|\sigma\|$ be as in Lemma~\ref{lemma_CG}, and $A$ be a symmetric positive definite matrix with spectrum $\lambda_1 < \lambda_2 < .... \leq b$.
  \begin{description}
\item[(a)] If the integers $n\geq 2$ and $d=d_n \in \{ 0,1,...,n-2 \}$ are such that
    \begin{eqnarray}
    &&\label{eq1.Th_super}
           \mbox{for $j=1,2,...,d$:} \quad \sigma((-\infty,\lambda_j])\geq j/N,
    \\ &&\label{eq2.Th_super}
           \lambda_{d}< a(n/N) \leq \lambda_{d+1}  \quad \mbox{(or $a(n/N)\leq \lambda_1$ in case $d=0$),}
    \end{eqnarray}
    then
    $$
       E_n(\Lambda(A)) \leq E_{n}([\lambda_{d+1},b])
       \leq
         \exp\Bigl(C_{BW} - N \int_0^{n/N} g_{S(t)}(0,\infty) dt \Bigr),
  $$
  and thus Conjecture~\ref{conjecture_superlinear} holds.
  \item[(b)] The above choice \eqref{eq2.Th_super} of $d$ is nearly optimal in the following sense:
      consider diagonal $A_N$ with eigenvalues satisfying $\sigma((-\infty,\lambda_{j,N}])= j/N$ for $j=1,...,N$.
      Furthermore, let $d=d_{n,N}$ with $\lambda_{d,N}< a(n/N) \leq \lambda_{d+1,N}$,
      then\footnote{We write $E_{n,N}$ instead of $E_n$ in order to indicate that here we consider the spectrum of $A_N$ depending on $N$.}
      $$
           \lim_{{n,N\to \infty}\atop {n/N\to t}} E_{n,N}([\lambda_{d_{n,N}+1,N},b])^{1/n}
           =
           \liminf_{{n,N\to \infty}\atop {n/N\to t \atop 0\leq d <n-1}} E_{n,N}([\lambda_{d+1,N},b])^{1/n}.
      $$
\end{description}
\end{Th}

It is also interesting to compare Theorem~\ref{Th_super}(a),(b) with  \cite[Theorem~3.1]{Beck-sharpness} which showed under the sole assumption \eqref{weak_limit}
(and for quite general measures $\sigma$) that, for any fixed compact set $S$, the quantity $E_{n,N}(S)^{1/n}$ is asymptotically greater than or equal
to the right-hand side of \eqref{eq.BB_ABK}. One of the consequences of our Theorem~\ref{Th_super} is that, roughly,  we can achieve equality for the interval $S=S(n/N)$.

Our proof of Theorem~\ref{Th_super} will be presented in \S\ref{sec_proof_th_super}, let study here some examples.

\begin{Exa}\label{example11}
   Consider the probability density
   $$
         \frac{d\sigma}{dx}(x) = \frac{1}{2\sqrt{1-x}} \quad\mbox{on} \quad
         [a(0),b]=[0,1] .
   $$
   For this measure we may apply Lemma~\ref{lemma_CG}, and a small computation shows for $0<t<\|\sigma\|=1$ that $a(t)=t^2$. We may also compute eigenvalues
   $\lambda_j$ satisfying equality in \eqref{eq1.Th_super}:
   $$
        \sigma([0,\lambda_j])=\frac{j}{N}
        \quad \mbox{iff}
        \quad \lambda_j = \frac{j}{N} \Bigl( 2 - \frac{j}{N}\Bigr),
   $$
   which behave like equidistant points for $j \ll N$.
   These are the eigenvalues used in Figure~\ref{fig1} and Figure~\ref{fig2}.
   In this special example we even have an explicit formula for the quantity $d=d_n$ of Theorem~\ref{Th_super}(a), namely
   $$
        d_n+1 = \lceil N ( 1 - \sqrt{1-(n/N)^2} )\rceil \approx \lceil \frac{n^2}{2N} \rceil,
   $$
   in particular $d_n=0$ for $n\leq 45$, $d_n=1$ for $46 \leq n \leq 64$, and $d_n=2$ for $65 \leq n \leq 78$,
   in accordance with the right-hand plot of Figure~\ref{fig2}.
\end{Exa}

In the previous example the small eigenvalues were approximately equidistant, with stepsize $2/N$, and the convex hull of the spectrum given approximately by $[2/N,1]$.
Up to correct scaling, a similar behavior is true for the eigenvalues of the finite difference discretization of the 2D Laplacian on the unit square
with Dirichlet boundary conditions, and thus
 the convergence curves should be similar. However, this is no longer true for higher dimensions $D \geq 3$,
where we expect that $\sigma'(x)$ grows like a constant times $x^{(D-2)/2}$ for small $x$, which motivates the following example.

\begin{figure}[t]
   \centerline{\includegraphics[width=0.49\textwidth]{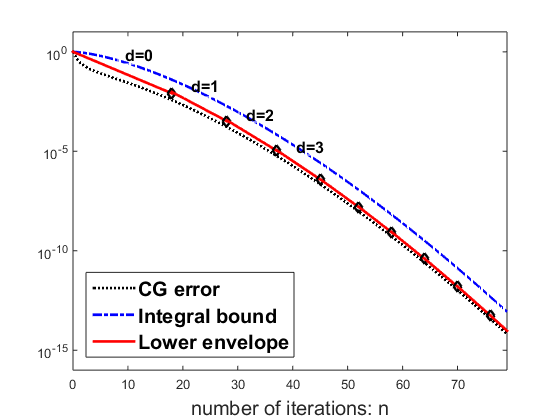}~
   \includegraphics[width=0.49\textwidth]{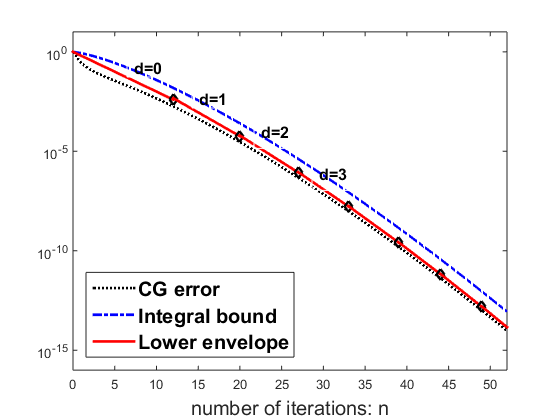}}
   \caption{Illustration of Theorem~\ref{Th_super}(a), with $\lambda_j$ for $j=1,...,N=1000$ as in Example~\ref{example12},
	 where on the left $\beta=0.5$ and on the right $\beta=1$.
   As lower bound we draw on both plots the relative CG error in energy norm, with $A,c,x_0^{CG}$ as in Figure~\ref{fig1} (blue dotted line).
	 The polygons are obtained from Theorem~\ref{Th_super}(a) (red solid line),
	 where we indicate in the plot the correspondence between a segment and the choice of $d$. We also draw the integral bound (blue dash-dotted line),
	showing that Conjecture~\ref{conjecture_superlinear} holds with $C=0$.
	 }\label{fig3}
\end{figure}

\begin{Exa}\label{example12}
   For a parameter $\gamma>0$, consider the density
   $$
       \frac{d\sigma}{dx}(x) = \frac{\gamma x^{\beta}}{\sqrt{b-x}} \quad\mbox{on} \quad
         [0,b] .
   $$
   In this example we only consider probability measures $\sigma$ and thus $\gamma b^{\beta+1/2} B(\beta+1,1/2)=1$, with $B(\cdot,\cdot)$ the beta function.
	 Notice that, for $\beta=0$, we recover Example~\ref{example11}.
   A small computation using Lemma~\ref{lemma_CG} gives
	 $$
	     a(t)/b = t^{\frac{1}{\beta+1/2}}.
	 $$
   We again choose $\lambda_j$ for $j=1,2,...,N=1000$ attaining equality in \eqref{eq1.Th_super}, however, there are no longer explicit formulas,
	 and thus the $\lambda_j$ have to be computed numerically.
   In Figure~\ref{fig3} we have plotted two examples for $b=1$, on the left for $\beta=0.5$ and on the right for $\beta=1$,
	 where in both cases we have chosen the approximately optimal $d=d_n$ of Theorem~\ref{Th_super}(a), in accordance with the statement of Theorem~\ref{Th_super}.
	 Notice also the well-known phenomena that the convergence of CG improves dramatically with $\beta$ getting larger.
\end{Exa}

\subsection{Structure of the paper}

The reminding of the paper is organized as follows.
Section 2 contains our results on discretizing the logarithmic potential of a class of measures including the extremal measure of Theorem~\ref{th_main}.
We first state our main Theorem~\ref{theorem-potential}, and then report in \S\ref{sec_distretize} about related results of Totik and of Lubinsky,
and about the link with weighted quadrature formulas. Subsequently, we give in \S\ref{proof_potential}
the structure of the proof of Theorem~\ref{theorem-potential}, where following Totik we write the discretization error as a sum of three sums.
We then state our original approach for dealing with these three sums, namely the mean value property of Theorem~\ref{Theorem:mean:property},
and describe in \S\ref{proof_bounding_sums} how to bound each of the three sums, with explicit constants.

In the third section we explain how to deduce Theorem~\ref{th_main} from Theorem~\ref{theorem-potential} (\S\ref{subsec_proof_th_main}),
and Theorem~\ref{Th_super} from Theorem~\ref{th_main} (\S\ref{sec_proof_th_super}). Subsequently, we give some concluding remarks.
Our (quite technical) proof of Theorem~\ref{Theorem:mean:property} is postphoned to Appendix~\ref{annexeA},
and in Appendix~\ref{annexeB} we gather some further technical results for dealing with our three sums.


\section{Discretization of a potential}\label{sec_discretization}

Our proof of Theorem~\ref{th_main} is based on the approximation of $kU^\mu$ with $\mu$ the equilibrium measure as in Theorem~\ref{th_main}
by $- \log | P_k(z)|$ with $P_k$ a suitable monic polynomial of degree $k$. We will show the following.

\begin{Th}
  \label{theorem-potential}
  Consider a measure $\mu\in \mathcal M_1([a,b])$ which has the density $$
       k\frac{d\mu}{dx}(t)= g(t)\frac{k}{\pi \sqrt{(t-a)(b-t)}}
  $$ for a function $g$ which is non negative, concave and
  increasing\footnote{
In particular, $g$ is continuous  and bounded on $(a,b)$, thus we may extend $g$ to become a continuous, non-negative, concave and increasing function in $[a,b]$.
}
  on $(a,b)$,
  such that $t \mapsto \frac{g(t)}{t-a}$ is convex on $(a,b)$.
	Then there exists a universal explicit constant $C_{BW}$ such that
	for each $k \geq 2$
    we may construct a monic polynomial $P_k$ of degree $k$ such that
  \begin{enumerate}
	  \item[(a)] $\forall z \in \mathbb C$:
     $\log |P_k(z)| + k
          U^{\mu}(z) \leq C_{BW}$,
    \item[(b)] $\forall x \in \mathbb R \setminus (a,b):$ $\log |P_k(x)| + k U^{\mu}(x)  \geq 0$.
  \end{enumerate}
\end{Th}

We will show in the proof of Theorem~\ref{th_main} that the extremal measure $\mu$ of Theorem~\ref{th_main} satisfies the assumptions of Theorem~\ref{theorem-potential}.
\begin{Exa}\label{example3}
   Another class of functions $g$ satisfying the assumptions of Theorem~\ref{theorem-potential} for $[a,b]=[-1,1]$ is given by
   \begin{eqnarray*}
          g(x) &=& (x+1)^\theta \pi / \int_{-1}^1 (t+1)^{\theta-1/2} (1-t)^{-1/2} \mathrm{d}t
          \\&=&
          \frac{\pi}{2^\theta} \frac{\Gamma(\theta+1)}{\Gamma(1/2) \Gamma(\theta+1/2) } (x+1)^\theta
   \end{eqnarray*}
   for $\theta \in [0,1]$.
\end{Exa}

We will describe in \S 2.1 related work for discretizing potentials under various assumptions,
but here the constants in general depend on $\mu$, see for instance \cite[\S VI.4]{Saff-Totik} for a summary.
In \S 2.2 we give a proof of Theorem~\ref{theorem-potential},
where we initially follow the approach of Totik in \cite[\S 2 and \S 9]{Totik:LNM},
see also the very accessible reference \cite[Method 1]{Lubinsky-Totik} for the particular case $g(t)=2t$ on $[a,b]=[0,1]$ (up to a quadratic change of variables).
Subsequently, we give in \S 2.3 a proof of three upper bounds we used in \S 2.2.
Since the general case follows from a linear change of variables, we will suppose in what follows that $[a,b]=[-1,1]$ in Theorem~\ref{theorem-potential}.

\subsection{How to discretize a potential?}\label{sec_distretize}

It is natural to approach the logarithmic potential $U^\mu(x)=\int \log(1/|x-t|) \mathrm{d}\mu(t)$ by a quadrature rule of the form
\begin{equation} \label{midpoint0}
     \sum_{j=0}^{k-1} \log\frac{1}{|x-\xi_j|} = - \log| P_k(x)|, \quad
     P_k(x)= \prod_{j=0}^{k-1} (x-\xi_j) ,
\end{equation}
for instance a weighted rectangular or midpoint rule, where we first cut $[-1,1]$ into $k$ subintervals $[t_j,t_{j+1}]$, $-1=t_0<t_1<...<t_k=1$,
of equal mass $\mu([t_j,t_{j+1}])=1/k$, and chose $\xi_j \in [t_j,t_{j+1}]$ for $j=0,...,k-1$.
As long as $x\not\in [-1,1]$ and the density of $\mu$ does not vary too much, we may bound the error $kU^\mu(x)+\log|P_k(x)|$ above and below,
and may even show convergence to $0$ for $k\to \infty$ for suitable choices of $\xi_j$.
In our case we have the additional difficulties that the density of $\mu$ may have singularities at $\pm 1$,
showing that the interval lengths $t_{j+1}-t_j$ may strongly vary in size for $j=0,1,...,k-1$,
and in addition in case $x\in [-1,1]$ we have to deal with a logarithmic singularity of the integrand.

Totik in \cite[Method 1]{Lubinsky-Totik} used the weighted midpoint rule
\begin{equation} \label{midpoint}
       \xi_j = \int_{t_j}^{t_{j+1}} t \mathrm{d}\mu(t) / \int_{t_j}^{t_{j+1}} \mathrm{d}\mu(t)
       = k \int_{t_j}^{t_{j+1}} t \mathrm{d}\mu(t)
\end{equation}
for $j=0,1,...,k-1$.
In the particular case $[a,b]=[0,1]$ and $g(t)=2t$, a proof of Theorem~\ref{theorem-potential} can be found in \cite[\S 2]{Lubinsky-Totik},
which strongly relies on the explicit knowledge of asymptotics for the points $\xi_j$ and $t_j$ as a function of $j$ and $k$ for $k \to \infty$,
and thus on the explicit knowledge of $\mu$. In \cite[Theorem~VI.4.2]{Saff-Totik} (see also the related result \cite[Lemma~9.1]{Totik:LNM}
where the roots of $P_k$ are slightly shifted into the complex plane), Totik considered probability measures $\mu$ with densities which are
continuous up to a finite number of singularities of the form $|t-a_j|^{\delta_j}$ for $\delta_j>-1$. These assumptions are true in the setting of Theorem~\ref{theorem-potential}.
He then shows the existence of (non explicit) constants $C_{T,1},C_{T,2}$
depending on $\mu$ but not on $k$ such that, for all $x\in \mathbb R$,
\begin{eqnarray*} &&
     \log |P_k(x)| + k
          U^{\mu}(x) \leq C_{T,1},
    \\&&
    \log |P_k(x)| + k U^{\mu}(x)  \geq C_{T,2}
    + \max \Bigl\{ 0 , \log(\mbox{dist}(x,\{ \xi_0,...,\xi_{k-1} \})) \Bigr\}.
\end{eqnarray*}
We see that the first inequality is as in Theorem~\ref{theorem-potential}(a), whereas the second one is clearly weaker than
Theorem~\ref{theorem-potential}(b) for $x\in \mathbb R \setminus (-1,1)$ close to $[-1,1]$,
since we get an additional term $\log(1/k^\beta)$ for some $\beta>0$.
Again, a proof of these statements uses heavily asymptotics for the points $\xi_j$ and $t_j$ as a function of $j$ and $k$ for $k \to \infty$,
and thus quite a bit of information on $\mu$.

Another technique of discretization has been considered by Lubinsky \& Levin
in \cite{Levin-Lubinsky1} and \cite{Levin-Lubinsky2}, see also the very accessible reference \cite[Method 2]{Lubinsky-Totik}
for the particular case $g(t)=2t$ on $[a,b]=[0,1]$ (up to a quadratic change of variables). With $t_0,...,t_k$ as before,
consider intermediate abscissa $t_{j+1/2}\in (t_j,t_{j+1})$ such that all intervals $[t_{j/2},t_{(j+1)/2}]$ have the same mass $1/(2k)$.
Given $x_0\in \mathbb R$, the authors then apply trapezian rule on most of the subintervals $[t_{j-1/2},t_{j+1/2}]$
corrected with suitable rectangle rules on the remaining 2 or 3 subintervals such that
$\{ \xi_0,...,\xi_{k-1} \} \subset \{ \pm 1,t_{1/2},t_{3/2},...,t_{k-1/2} \}$. Up to a (quadratic) change of variables,
the authors of \cite[Theorem 9.1]{Levin-Lubinsky1} suppose that
$$
    \frac{d\mu}{dx}(t)=\frac{(t+1)h(t)}{\pi \sqrt{1-t^2}}
$$ with $h$ continuous and $>0$ on $[-1,1]$, and the modulus of continuity satisfies that $\log(1/\delta)\omega(h,\delta)$ is bounded above by some
$\Gamma>0$ for $\delta\in (0,1)$. In this case, for all $x\in \mathbb R$,
\begin{eqnarray*} &&
     \log |P_k(x)| + k
          U^{\mu}(x) \leq C_{LL,1},
    \\&&
    \log |P_k(x_0)| + k U^{\mu}(x_0)  \geq C_{LL,2},
\end{eqnarray*}
where $C_{LL,1},C_{LL,2}$ are (non explicit) constants depending only on $\Gamma$ and the minimum and maximum of $h$ on $[-1,1]$.
Note that the assumptions of \cite[Theorem 9.1]{Levin-Lubinsky1} and those of Theorem~\ref{theorem-potential} are different and do not imply each other,
see for instance Example~\ref{example3} for $\theta<1$. However, the above inequalities are quite close to those of Theorem~\ref{theorem-potential},
though our constants do not depend on $\mu$, and our $P_k$ does not depend on $x_0$, and we only allow $x_0\in \mathbb R \setminus (a,b)$.

\begin{Exa}\label{example4}
   In the particular case $[a,b]=[-1,1]$ and $g=1$ in Theorem~\ref{theorem-potential}, we have explicit formulas
   $$
        t_j = - \cos(\pi \frac{j}{k}) , \quad
        \xi_j = -c_k \cos(\pi \frac{2j+1}{2k}) , \quad c_k=\frac{2k}{\pi} \sin(\frac{\pi}{2k}).
   $$
   Here the midpoint approach of Totik gives the monic polynomial
   $$
       P_k(x) =  2 ( \frac{c_k}{2} )^k T_k(x/c_k)
   $$
   which is not optimal for the one-sided approximation of $kU^\mu(x)$ in Theorem~\ref{theorem-potential} or the sharpness of the classical
	Bernstein-Walsh inequality as discussed in Example~\ref{example1}, but good enough for concluding in Theorem~\ref{theorem-potential}.
\end{Exa}

The previous example is misleading in the sense that in general there is no such sufficiently explicit formula for the $t_j$ nor
the $\xi_j$ which will allow us to conclude in Theorem~\ref{theorem-potential}.

\subsection{Structure of the proof of Theorem~\ref{theorem-potential}}\label{proof_potential}

We start by observing that, with the choices \eqref{midpoint0}, \eqref{midpoint},
$$
     \log |P_k(x)| + k
          U^{\mu}(x) =  k \sum_{j=0}^{k-1} \int_{t_j}^{t_{j+1}}
          \log \left|\frac{x-\xi_j}{x-t}\right| \, \mathrm{d}\mu(t).
$$
The following classical lemma shows Theorem~\ref{theorem-potential}(b).

\begin{Lemma}\label{lem_quadrature}
    $$
          k \int_{t_j}^{t_{j+1}}
          \log \left| \frac{x-\xi_j}{x-t} \right| \, \mathrm{d}\mu(t)
          \left\{\begin{array}{ll}
            \geq 0 & \mbox{for $x\in \mathbb R \setminus (t_j,t_{j+1})$},
            \\
            \displaystyle
           \leq \frac{1}{4} \frac{(t_{j+1}-t_j)^2}{(x-t_j)(x-t_{j+1})}
           & \mbox{for $x\in \mathbb R \setminus [t_j,t_{j+1}]$}.
        \end{array} \right.
    $$
\end{Lemma}
\dem{
    Using the fact that $m(t)=\log |\frac{x-\xi_j}{x-t}|$ is convex on $[t_j,t_{j+1}]$ by assumption on $x$, we know
    that $m(t)\geq m(\xi_j)+m'(\xi_j)(t-\xi_j)
    = m'(\xi_j)(t-\xi_j)$, and thus
    $$
         k \int_{t_j}^{t_{j+1}}
          \log \left| \frac{x-\xi_j}{x-t} \right| \, \mathrm{d}\mu(t)
          \geq m'(\xi_j) \int_{t_j}^{t_{j+1}}
          (t-\xi_j) \, \mathrm{d}\mu(t) = 0,
    $$
    where in the last equality we have used \eqref{midpoint}.
    Also, using the convexity of $m$ and the inequality $\log(x) \leq x-1$ we obtain
    \begin{eqnarray*}
         m(t) & \leq & m(t_j) \frac{t_{j+1}-t}{t_{j+1}-t_j} + m(t_{j+1}) \frac{t-t_j}{t_{j+1}-t_j} \\
         & \leq & \frac{t_j-\xi_j}{x-t_j} \frac{t_{j+1}-t}{t_{j+1}-t_j} + \frac{t_{j+1}-\xi_j}{x-t_{j+1}} \frac{t-t_j}{t_{j+1}-t_j}.
    \end{eqnarray*}
    Integrating and using again \eqref{midpoint} we conclude that
    \begin{eqnarray*}
         k \int_{t_j}^{t_{j+1}} \log | \frac{x-\xi_j}{x-t} | \, \mathrm{d}\mu(t)
         & \leq & \frac{t_j-\xi_j}{x-t_j} \frac{t_{j+1}-\xi_j}{t_{j+1}-t_j} + \frac{t_{j+1}-\xi_j}{x-t_{j+1}} \frac{\xi_j-t_j}{t_{j+1}-t_j} \\
         & = &
         \frac{(t_{j+1}-\xi_j)(\xi_j-t_j)}{(x-t_j)(x-t_{j+1})}
         \leq \frac{1}{4}
         \frac{(t_{j+1}-t_j)^2}{(x-t_j)(x-t_{j+1})}.
    \end{eqnarray*}
}

\begin{Rem}\label{remark_sharpness}
   (a) The interested reader might have noticed that, by the same argument, the inequality of Theorem~\ref{theorem-potential}(b), namely $\log |P_k(x)| + k
          U^{\mu}(x)\geq 0$, also holds for $x\in \{ t_0,t_1,...,t_k \}$.
   \\ (b) For $x>1$ (and similarly for $x<-1$), the right-hand side of Theorem~\ref{theorem-potential}(b) cannot be improved since, by Lemma~\ref{lem_quadrature} and Lemma~\ref{lemma:general:case:racine:t1-t0}(c),
   \begin{eqnarray*}
        \log |P_k(x)| + k
          U^{\mu}(x)
          &\leq& \max_{\ell=0,...,k-1} \frac{t_{\ell+1}-t_\ell}{4}
          \sum_{j=0}^{k-1}\frac{t_{j+1}-t_j}{(x-t_j)(x-t_{j+1})}
         \\&=& \max_{\ell=0,...,k-1} \frac{t_{\ell+1}-t_\ell}{2(x^2-1)}
         \leq \frac{1}{(x^2-1)} \left( \frac{3\pi}{2k}\right)^{1/3} .
   \end{eqnarray*}
   (c) For $x\in \mathbb C \setminus \mathbb R$, it is not too difficult to show that $m(t)=\log |\frac{x-\xi_j}{x-t}|$ satisfies
   $$
        | m(t) - m(\xi_j) - (t-\xi_j) m'(\xi_j) | \leq \frac{(t_{j+1}-t_j)^2}{2 \, \mbox{dist}(x,[-1,1])^2} ,
   $$
   and hence by Lemma~\ref{lemma:general:case:racine:t1-t0}(c)
   \[
       | \log |P_k(x)| + k
          U^{\mu}(x)
       | \leq \sum_{j=0}^{k-1} \frac{(t_{j+1}-t_j)^2}{2 \, \mbox{dist}(x,[-1,1])^2}
       \leq \frac{1}{\mbox{dist}(x,[-1,1])^2} \left( \frac{12\pi}{k}\right)^{1/3}.
   \]
   Thus, for sufficiently large $k$, the inequality of Theorem~\ref{theorem-potential}(b) also holds for non-real $x$ up to some arbitrarily small constant.
\end{Rem}

Let us now turn to a proof of Theorem~\ref{theorem-potential}(a).
We claim that it is sufficient to show Theorem~\ref{theorem-potential}(a) for $x\in [-1,1]=\supp(\mu)$, since then for $\mu$-almost all $x$
$$
     kU^\mu(x)  \leq C_{BW} - \log| P_k(z)| = C_{BW} + \sum_{j=0}^{k-1} U^{\delta_{\xi_j}}(x),
$$
and thus this inequality holds for all $x\in \mathbb C$ by the principle of domination \cite[Theorem~II.3.2]{Saff-Totik} and the finiteness of $I(\mu)$.
Therefore, let $x\in [-1,1]$ and, more precisely,
\begin{equation} \label{eq.choice}
        j_0\in \{ 0,1,...,k-1\} \quad \mbox{with} \quad x\in [t_{j_0},t_{j_0+1}].
\end{equation}
According to Lemma~\ref{lem_quadrature}, we get the following upper bound
\begin{equation} \label{eq.estimate}
             \log |P_k(x)| + k U^{\mu}(x) \leq \Sigma_1 + \Sigma_2 + \Sigma_3
\end{equation}
with
\begin{align*}
  {\sum}_1 & = \sum_{j=0}^{j_0-2} \int_{t_{j}}^{t_{j+1}} \log \left| \frac{x-\xi_j}{x-t} \right| k \, \mathrm{d}\mu(t)
  \leq \frac{1}{4} \sum_{j=0}^{j_0-2} \frac{(t_{j+1}-t_j)^2}{(t_{j_0}-t_{j+1})^2}, \\
  {\sum}_2 & = \sum_{j=\max\{0,j_0-1\}}^{\min\{j_0+1,k-1\}} \int_{t_j}^{t_{j+1}} \log \left| \frac{x-\xi_j}{x-t} \right| k \, \mathrm{d}\mu(t),  \\
  {\sum}_3 & = \sum_{j=j_0+2}^{k-1} \int_{t_{j}}^{t_{j+1}} \log \left| \frac{x-\xi_j}{x-t} \right|  k \, \mathrm{d}\mu(t)
  \leq \frac{1}{4} \sum_{j=j_0+2}^{k-1} \frac{(t_{j+1}-t_j)^2}{(t_{j}-t_{j_0+1})^2} .
\end{align*}

Already in the particular Chebyshev case of Example~\ref{example4} one may check that such a simple telescop sum trick as in Remark~\ref{remark_sharpness}
does not allow to conclude, since in general $|t_j-t_\ell|$ does not behave uniformly for $j,\ell \in \{ 0,1,...,k-1 \}$ like $|j-\ell|/k$,
as it would be the case for equidistant points. We will discuss our upper bounds for the above three sums
in the Propositions~\ref{general:case:sum:3}--\ref{general:case:sum:2} of \S \ref{proof_bounding_sums},
which allows us to conclude the proof of Theorem~\ref{theorem-potential}, with the explicit constant
$$
    C_{BW} = 872 + 32 + 14 479 = 15 383.
$$
So far we followed quite closely the reasoning in the literature, with more explicit constants.
In all considerations to follow we will require precise lower and upper bounds for the ratio
$$
       \frac{j-\ell}{t_j-t_\ell}
$$
which will follow from a new mean value property for the cumulative distribution function
\begin{equation} \label{distribution}
      W_g(x)= k \int_{-1}^x \mathrm{d}\mu(t), \quad W'_g(t) = g(t) W_1'(t) = \frac{k g(t)}{\pi \sqrt{1-t^2}},
\end{equation}
since $W_g(t_j)=j$ for $j=0,1,...,k$.

\begin{Th}
  \label{Theorem:mean:property}
  Under the assumptions of Theorem~\ref{theorem-potential} with $[a,b]=[-1,1]$, we have for all distinct $x,t\in [-1,1]$
  \begin{equation}
     c_1 W_g'(\frac{t+x}{2}) \leq \frac{W_g(t)-W_g(x)}{t-x}  \leq c_2 W_g'(\frac{t+x}{2}) .
  \end{equation}
  where $c_1 = \frac{1}{4}$ and $c_2 = \pi\sqrt{2}$.
\end{Th}

Notice that, even for the particular case $g=1$ and $W_1(-\cos(\alpha)) = k\alpha/\pi$, this statement is not totally obvious, but can be verified
by means of elementary computations with improved constants $c_1$ and $c_2$, see Lemma~\ref{lemma:mean:property:Chebyshev:case} below.
The proof for general $g$ is strongly based on Jensen's inequality, we refer the reader to Appendix~\ref{annexeA} for details.

\subsection{Bounding three sums}\label{proof_bounding_sums}

For concluding our proof of Theorem~\ref{theorem-potential}, it remains to obtain upper bounds for the three terms on the right-hand side of \eqref{eq.estimate}, where we will proceed in order of difficulty, and apply beside Theorem~\ref{Theorem:mean:property} a certain number of technical results established in Appendix~\ref{annexeB}, and recalled below. In the reminder of this section we will always suppose that the assumptions of Theorem~\ref{theorem-potential} hold with $[a,b]=[-1,1]$ and $j_0$ is chosen as in \eqref{eq.choice}.

We start with the sum
\[ {\sum}_3 \leq \frac{1}{4} \sum_{j=j_0+2}^{k-1} \frac{(t_{j+1}-t_j)^2}{(t_{j}-t_{j_0+1})^2} , \]
where beside Theorem~\ref{Theorem:mean:property} we rely on an upper bound for the quantity
\[ (1+\frac{t_j + t_{j+1}}{2})\Bigl/(1+t_j) , \]
see Lemma \ref{lemma:rapport:I_j:avant}.
\begin{Prop}
  \label{general:case:sum:3}
  There holds
  \begin{equation*}
    {\sum}_3 \leq  \frac{c_2^2 c_5\pi^2}{12c_1^2} \leq 872 .
  \end{equation*}
\end{Prop}
\dem{
By Theorem~\ref{Theorem:mean:property}
\[ {\sum}_3 \leq \frac{c_2^2}{4c_1^2} \sum_{j=j_0+2}^{k-1}  \frac{1}{(j-j_0-1)^2} \frac{W_g'(\frac{t_j+t_{j_0+1}}{2})^2}{W_g'(\frac{t_j+t_{j+1}}{2})^2} . \]
As $g$ is increasing and $j>j_0+1$, we have that
$g(\frac{t_j+t_{j_0+1}}{2})\leq g(\frac{t_j+t_{j+1}}{2})$, and thus
\begin{align*}
  \frac{W_g'(\frac{t_j+t_{j_0+1}}{2})^2}{W_g'(\frac{t_j+t_{j+1}}{2})^2}
  & = \frac{g(\frac{t_j+t_{j_0+1}}{2})^2}{g(\frac{t_j+t_{j+1}}{2})^2}  \frac{W'_1(\frac{t_j+t_{j_0+1}}{2})^2}{W'_1(\frac{t_j+t_{j+1}}{2})^2} \\
  & \leq \frac{W_1'(\frac{t_j+t_{j_0+1}}{2})^2}{W_1'(\frac{t_j+t_{j+1}}{2})^2} \leq \frac{1+\frac{t_j+t_{j+1}}{2}}{1+\frac{t_j+t_{j_0+1}}{2}}
  \\
  & \leq 2 \frac{1+\frac{t_j+t_{j+1}}{2}}{1+t_j} \leq 2 c_5,
\end{align*}
where in the last inequality we have applied Lemma~\ref{lemma:rapport:I_j:avant}. Combining these two results yields the claimed upper bound.
}

Let us now turn to the sum
\[ {\sum}_1 \leq \frac{1}{4} \sum_{j=0}^{j_0-2} \frac{(t_{j+1}-t_j)^2}{(t_{j_0}-t_{j+1})^2}  . \]
Here we require beside Theorem \ref{Theorem:mean:property} also upper bounds for the two ratios
$$
      \frac{1+\frac{t_j + t_{j+1}}{2}}{1+t_{j+1}} , \quad
      \mbox{and} \quad
      \frac{(j+1)^2}{j_0^2} \frac{1+\frac{t_{j+1}+t_{j_0}}{2}}{1+t_{j+1}}
$$
for $j \leq j_0-1\leq k-2$, see Lemma~\ref{lemma:rapport:I_j:apres} and
Lemma~\ref{lemma:rapport:I_j:dur}.

\begin{Prop}
  \label{general:case:sum:1}
  There holds
  \begin{equation}
    {\sum}_{1} \leq \frac{c_2^2\pi^2}{6c_1^2} ( 18+\pi^2 ) \leq 14 479.
  \end{equation}
\end{Prop}
\dem{
  Using the fact that $j<j_0-1$, and that $g(t)=(1+t)h(t)$ with a decreasing function $h$, we find that
  \begin{align*}
   \frac{g(\frac{t_{j+1}+t_{j_0}}{2})}{g(\frac{t_{j}+t_{j+1}}{2})} & = \frac{1+\frac{t_{j+1}+t_{j_0}}{2}}{1+\frac{t_j+t_{j+1}}{2}} \frac{h(\frac{t_{j+1}+t_{j_0}}{2})}{h(\frac{t_{j}+t_{j+1}}{2})}  \\
   & \leq \frac{1+\frac{t_{j+1}+t_{j_0}}{2}}{1+\frac{t_j+t_{j+1}}{2}}.
  \end{align*}
  This allows us to write
  \begin{align*}
   \frac{W_g'(\frac{t_{j+1}+t_{j_0}}{2})^2}{W_g'(\frac{t_j+t_{j+1}}{2})^2} & \leq \frac{(1+\frac{t_{j+1}+t_{j_0}}{2})^2}{(1+\frac{t_j+t_{j+1}}{2})^2} \frac{W_1'(\frac{t_{j+1}+t_{j_0}}{2})^2}{W'_1(\frac{t_j+t_{j+1}}{2})^2}  \\
   & \leq \frac{1+\frac{t_{j+1}+t_{j_0}}{2}}{1+\frac{t_j+t_{j+1}}{2}} \frac{1-\frac{t_j+t_{j+1}}{2}}{1-\frac{t_{j+1}+t_{j_0}}{2}} \\
   & \leq 2 \frac{1+\frac{t_{j+1}+t_{j_0}}{2}}{1+t_{j+1}} \frac{1-\frac{t_j+t_{j+1}}{2}}{1-\frac{t_{j+1}+t_{j_0}}{2}} ,
\end{align*}
  and thus, again by Theorem \ref{Theorem:mean:property},
  \begin{align*}
  {\sum}_1 & \leq \frac{c_2^2}{2c_1^2} \sum_{j=0}^{j_0-2} \frac{1}{(j_0-j-1)^2} \frac{1+\frac{t_{j+1}+t_{j_0}}{2}}{1+t_{j+1}} \frac{1-\frac{t_j+t_{j+1}}{2}}{1-\frac{t_{j+1}+t_{j_0}}{2}}.
  \end{align*}
  The following arguments depend on the sign of $t_{j+1}$.
  We therefore set $j_1=j_0-1$ if $t_{j_0-1}<0$, and else chose $j_1\in \{ 0,1,...,j_0-2\}$ with $t_{j_1}<0\leq t_{j_1+1}$,
  and cut our sum into two parts $\Sigma_1=\Sigma_{1,1}+\Sigma_{1,2}$,
  where in the first sum $j\in \{ j_1,...,j_0-2\}$, and in the second one $j\in \{0,...,j_1-1\}$.

  If $j\geq j_1$ and thus $t_{j+1}\geq 0$,
  \begin{align*}
        \frac{1+\frac{t_{j+1}+t_{j_0}}{2}}{1+t_{j+1}} \frac{1-\frac{t_j+t_{j+1}}{2}}{1-\frac{t_{j+1}+t_{j_0}}{2}}
        \leq 2 \frac{1-\frac{t_j+t_{j+1}}{2}}{1-\frac{t_{j+1}+t_{j_0}}{2}}
        & \leq 4 \frac{1-\frac{t_j+t_{j+1}}{2}}{1-t_{j+1}} \leq 36,
  \end{align*}
  where in the last inequality we have applied  Lemma~\ref{lemma:rapport:I_j:apres}. Hence,
  \begin{equation}
  \label{equation:sum:1:1}
  {\sum}_{1,1} \leq \frac{18c_2^2}{c_1^2} \sum_{j=j_1}^{j_0-2} \frac{1}{(j_0-j-1)^2} .
  \end{equation}

  If $j<j_1$ and thus $t_{j+1}< 0$,
    \begin{align*}
        \frac{1+\frac{t_{j+1}+t_{j_0}}{2}}{1+t_{j+1}} \frac{1-\frac{t_j+t_{j+1}}{2}}{1-\frac{t_{j+1}+t_{j_0}}{2}}
        &
        \leq
        4
                \frac{1+\frac{t_{j+1}+t_{j_0}}{2}}{1+t_{j+1}}
        \leq \frac{\pi^2}{2} \frac{j_0^2}{(j+1)^2},
  \end{align*}
  where the last inequality follows from Lemma~\ref{lemma:rapport:I_j:dur}. Thus,
  \begin{align*}
     {\sum}_{1,2}
     & \leq \frac{c_2^2\pi^2}{4c_1^2} \sum_{j=0}^{j_1-1} \frac{j_0^2}{(j_0-j-1)^2(j+1)^2} \\
     & \leq \frac{c_2^2\pi^2}{4c_1^2} \left(
     \sum_{j=0,j+1<j_0/2}^{j_1-1} \frac{4}{(j+1)^2}
     + \sum_{j=0,j+1\geq j_0/2}^{j_1-1} \frac{4}{(j_0-j-1)^2}
     \right).
  \end{align*}
  Since $\pi^2\leq 18$, a combination with \eqref{equation:sum:1:1} gives the upper bound for $\Sigma_1$ as claimed in Proposition~\ref{general:case:sum:1}.
}

We finally discuss in our third proposition the expression
\[
     {\sum}_2 = \sum_{j=\max\{0,j_0-1\}}^{\min\{j_0+1,k-1\}} \int_{t_j}^{t_{j+1}} \log \left| \frac{x-\xi_j}{x-t} \right| W_g'(t) \, \mathrm{d}t,
\]
where we integrate in a neighborhood of $x$ and thus have to deal with the logarithmic singularity of the integrand.
Here again Theorem~\ref{Theorem:mean:property} will be essential. As maybe expected from \cite{Totik:LNM},
our proof for $j_0\in \{ 1,2,...,k-2\}$ is quite different from that for $x$ close to the endpoints
and thus $j_0\in \{ 0,k-1\}$: in the first case, we require lower and upper bounds for the ratio of the lengths
of two consecutive intervals $[t_j,t_{j+1}]$ established in Lemma~\ref{lemma:rapport:I_j:I_j+1}, whereas in the second case we require upper bounds for
\[
   k \sqrt{t_1-t_0} ,\quad \mbox{and} \quad
   k \sqrt{t_{k}-t_{k-1}},
\]
see Lemma~\ref{lemma:general:case:racine:t1-t0}.

\begin{Prop}
  \label{general:case:sum:2}
  There holds
  \[
    {\sum}_{2} \leq 6 c_2 + \log \left( \frac{6c_2\sqrt{c_5}}{c_1} \right) \leq 32 .
  \]
\end{Prop}

\dem{
By integration by part,
\begin{align*}
  \int_{t_j}^{t_{j+1}} \log \left| \frac{x-\xi_j}{x-t} \right| W_g'(t) \mathrm{d}t = & \left[ \log \left| \frac{x-\xi_j}{x-t}  \right| \left( W_g(t)-W_g(x) \right) \right]_{t_j}^{t_{j+1}} \\
  & + \int_{t_j}^{t_{j+1}} \frac{W_g(t)-W_g(x)}{t-x}\mathrm{d}t .
\end{align*}
In order to make our formulas a bit easier to read, we write $j_1=\max\{ 0,j_0-1\}$, $j_2=\min\{k-1,j_0+1\}$, and get
$\Sigma_2=\Sigma_{2,1}+\Sigma_{2,2}$, with
\begin{align*}
   {\sum}_{2,1} &= \int_{t_{j_1}}^{t_{j_2+1}}
   \frac{W_g(t)-W_g(x)}{t-x}\mathrm{d}t ,
   \\
   {\sum}_{2,2} &= \sum_{j=j_1}^{j_2}
   \left(
   \log \left| \frac{x-\xi_j}{x-t_{j+1}}  \right| \left( W_g(t_{j+1})-W_g(x) \right)\right.
   \\&\quad \quad \quad +
   \left.
   \log \left| \frac{x-\xi_j}{x-t_{j}}  \right| \left( W_g(x)-W_g(t_{j}) \right)
   \right).
\end{align*}
   The first term is easily bounded. Indeed, using Theorem~\ref{Theorem:mean:property}, we get
   \begin{align*}
  {\sum}_{2,1} & \leq c_2 \int_{t_{j_1}}^{t_{t_{j_2+1}}} W_g'(\frac{t+x}{2}) \mathrm{d}t  \\
  & \leq 2c_2  \left( W_g(\frac{x+t_{j_2+1}}{2}) -W_g(\frac{x+t_{j_1}}{2}) \right) \\
  & \leq 2c_2  \left( W_g(t_{j_2+1}) -W_g(t_{j_1}) \right)
  = 2c_2 (j_2+1-j_1),
   \end{align*}
  and thus $\Sigma_{2,1} \leq 6 c_2$ for $j_0\in \{ 1,...,k-2\}$, and  $\Sigma_{2,1} \leq 4 c_2$ for $j_0\in \{ 0,k-1\}$.

  It remains to give an upper bound for $\Sigma_{2,2}$. We first study the case $j_0\in \{ 1,...,k-2\}$ and thus $j_1=j_0-1$, $j_2=j_0+1$.
  As $W_g(t_{j+1})-W_g(t_j)=1$ for every $j$, we notice that
   \begin{equation} \label{eq.rewriting}
W_g(t_{j_0+2})-W_g(x) = 2 (W_g(t_{j_0+1})-W_g(x)) + (W_g(x)-W_g(t_{j_0}))
  \end{equation}
  and
  \[ W_g(x)-W_g(t_{j_0-1}) = (W_g(t_{j_0+1})-W_g(x)) + 2 (W_g(x)-W_g(t_{j_0})) .\]
	Inserting this information into $\Sigma_{2,2}$, we obtain, after some elementary computations,
  \begin{align*}
  {\sum}_{2,2}
	= & \underbrace{\left[ W_g(t_{j_0+1}) - W_g(x) \right]}_{\geq 0} \log \frac{|x-\xi_{j_0}|}{|x-t_{j_0+2}|}
  \underbrace{\frac{|x-\xi_{j_0-1}|}{|x-t_{j_0-1}|}}_{\leq 1} \underbrace{\frac{|x-\xi_{j_0+1}|}{|x-t_{j_0+2}|}}_{\leq 1}  \\
  & + \underbrace{\left[ W_g(x) - W_g(t_{j_0}) \right]}_{\geq 0} \log \frac{|x-\xi_{j_0}|}{|x-t_{j_0-1}|}
  \underbrace{\frac{|x-\xi_{j_0-1}|}{|x-t_{j_0-1}|}}_{\leq 1} \underbrace{\frac{|x-\xi_{j_0+1}|}{|x-t_{j_0+2}|}}_{\leq 1} \\
	\leq & \left[ W_g(t_{j_0+1}) - W_g(x) \right] \log \left| \frac{x-\xi_{j_0}}{x-t_{j_0+2}} \right|
  + \left[ W_g(x) - W_g(t_{j_0}) \right] \log \left| \frac{x-\xi_{j_0}}{x-t_{j_0-1}} \right| \\
  \leq & \log \left( \frac{6c_2\sqrt{c_5}}{c_1} \right) .
  \end{align*}
  In order to justify the last inequality, we have to distinguish two cases.
  In case $x \in [\xi_{j_0},t_{j_0+1}]$, we find that
  $\log | \frac{x-\xi_{j_0}}{x-t_{j_0-1}}|\leq 0$, implying that
  \[
  {\sum}_{2,2} \leq \left[ W_g(x) - W_g(t_{j_0}) \right] \log \left| \frac{x-\xi_{j_0}}{x-t_{j_0-1}} \right| \leq \log \frac{t_{j_0+1}-t_{j_0}}{t_{j_0}-t_{j_0-1}},
  \]
  and we conclude with help of Lemma~\ref{lemma:rapport:I_j:I_j+1}.
  The case $x \in [t_{j_0},\xi_{j_0}]$ is similar, here ${\sum}_{2,2} \leq \log \frac{t_{j_0+1}-t_{j_0}}{t_{j_0+2}-t_{j_0+1}}$, and we conclude again using
  Lemma~\ref{lemma:rapport:I_j:I_j+1}.

Let us now consider the sum $\Sigma_{2,2}$ for the particular case $j_0=0$ and thus $j_1=0$, $j_2=1$. Using \eqref{eq.rewriting},
this sum can be bounded above as before by
\begin{align*}
  {\sum}_{2,2} = & \underbrace{\left[ W_g(t_{1}) - W_g(x) \right]}_{\geq 0} \log \frac{|x-\xi_{0}|}{|x-t_{2}|}
  \underbrace{\frac{|x-\xi_{1}|}{|x-t_{2}|}}_{\leq 1}
  \\
  & + \underbrace{\left[ W_g(x) - W_g(t_{0}) \right]}_{\geq 0} \log \frac{|x-\xi_{0}|}{|x-t_{0}|}
  \underbrace{\frac{|x-\xi_{1}|}{|x-t_{2}|}}_{\leq 1} \\
  \leq & \left[ W_g(t_{1}) - W_g(x) \right] \log \left| \frac{x-\xi_{0}}{x-t_{2}} \right|
  + \left[ W_g(x) - W_g(t_{0}) \right] \log \left| \frac{x-\xi_{0}}{x-t_{0}} \right| .
\end{align*}
We have to consider three different cases: if $x\in [\frac{t_{0}+\xi_{0}}{2},\frac{\xi_{0}+t_{1}}{2}]$ then $\Sigma_{2,2} \leq 0$.
If $x\in [\frac{\xi_{0}+t_{1}}{2},t_{1}]$, then
\[
   {\sum}_{2,2} \leq \left[ W_g(t_{1}) - W_g(x) \right] \log \left| \frac{x-\xi_{0}}{x-t_{2}} \right|\leq \log \frac{t_{1}-t_{0}}{t_{2}-t_{1}}
   \leq \log \left( \frac{6c_2\sqrt{c_5}}{c_1} \right)
\] as before. Finally, in the case $x\in [t_{0},\frac{t_{0}+\xi_{0}}{2}]$, we use the fact that $|x-\xi_{0}| \leq t_1-t_0$,
and apply Theorem~\ref{Theorem:mean:property} in order to get
\begin{align*}
  {\sum}_{2,2}   & \leq \left[ W_g(x) - W_g(t_{0}) \right] \log \left| \frac{t_1-t_0}{x-t_{0}} \right| \\
  & \leq c_2 (x-t_0) g(\frac{t_0+x}{2}) W_1'(\frac{t_0+x}{2}) \log \left| \frac{t_1-t_0}{x-t_{0}} \right| .
\end{align*}
Since $\frac{t_0+x}{2}\leq 0$ and  $\frac{t_0+x}{2}\leq t_1$, we have that
\[
   (x-t_0) g(\frac{t_0+x}{2}) W_1'(\frac{t_0+x}{2}) \leq \frac{k}{\pi} g(t_1) \frac{x-t_0}{\sqrt{1+\frac{t_0+x}{2}}} = \frac{k\sqrt{2}}{\pi} g(t_1)\sqrt{x-t_0} .
\]
Using the fact that $\max\limits_{y\geq 0} \sqrt{y} \log \frac{1}{y} = 2/e$,
we conclude with help of Lemma~\ref{lemma:general:case:racine:t1-t0}(a) that
\begin{align*}
  {\sum}_{2,2} & \leq \frac{2c_2 \sqrt{2}}{\pi e} k g(t_1) \sqrt{t_1-t_0}
  \leq  \frac{6}{e} c_2 .
\end{align*}
The reasoning for $j_0=k-1$ is similar and allows for the same conclusion, we just have to replace
Lemma~\ref{lemma:general:case:racine:t1-t0}(a) by Lemma~\ref{lemma:general:case:racine:t1-t0}(b) providing an upper bound for $k\sqrt{t_k-t_{k-1}}$. Thus
\[
     {\sum}_{2} = {\sum}_{2,1}+{\sum}_{2,2} \leq \max \left\{ 6c_2 + \log \left( \frac{6c_2\sqrt{c_5}}{c_1} \right), 4c_2 + \frac{6}{e} c_2  \right\} ,
\]
and the statement follows.
}

\section{Proof of the main theorems}\label{sec_proof_th_main}

\subsection{Proof of Theorem~\ref{th_main}}\label{subsec_proof_th_main}

Let us first show our claim \eqref{eq.extremal5} for the support of the equilibrium measure $\mu$.
We observe that the external field $Q=U^{\rho/k}$ is convex on $\Sigma=[\alpha,b]$ and hence $\supp(\mu)=[a,b']$ for some $\alpha \leq a < b'\leq b$ by \cite[Theorem~IV.1.10(b)]{Saff-Totik}. Since $U^\mu+Q$ is strictly decreasing on $(b',\infty)$, the equilibrium condition \eqref{eq.extremal2} tells us that necessarily $b=b'$.
We show below the two implications
\begin{eqnarray}
    \label{claim_alpha1}
     \mbox{for some $a>\alpha$}: && \mbox{$\supp(\mu)=[a,b]$  implies that $\eta(a)= k+\| \rho\|$,}
    \\&& \label{claim_alpha2}
     \mbox{$\supp(\mu)=[\alpha,b]$  implies that $\eta(\alpha)\leq k+\| \rho\|$,}
\end{eqnarray}
with the strictly decreasing $\eta$ as in \eqref{condition:mu:positive}.
Since there is exactly one solution $>\alpha$ of the equation $\eta(a)= k+\| \rho\|$ iff $\eta(\alpha) > k+\| \rho\|$, our statement on $\supp(\mu)$ follows.

For a proof of \eqref{claim_alpha1}, suppose that $\supp(\mu)=[a,b]$ for some $a>\alpha$.
Then, by \cite[Theorem~IV.1.11(ii)]{Saff-Totik},  the derivative of the $F$-functional of Mhaskar and Saff
$$
     y \mapsto \log\frac{b-y}{4} - \int Q d\omega_{[y,b]}
     = (1+\frac{\|\rho\|}{k})\log\frac{b-y}{4} + \frac{1}{k} \int g_{[y,b]}(x,\infty) d\rho(x)
$$
must vanish at $y=a$, and a small calculation gives
the necessary condition $$
      0 = \frac{1}{k(b-a)}\Bigl( k+\| \rho \| - \eta(a)\Bigr)
$$
and thus $\eta(a)=k+\| \rho \|$, implying \eqref{claim_alpha1}.

In order to show \eqref{claim_alpha2} together with the representation \eqref{eq.extremal5} of $f$, let $\supp(\mu)=[a,b]$ for some $a\in [\alpha,b)$.
We denote by $Bal(\rho,[a,b])$ the measure obtained by balayage onto the interval $[a,b]$, see \cite[\S II.4]{Saff-Totik}. Then, by construction,
$$
    k \mu + Bal(\rho,[a,b])
$$
is a positive measure of mass $k+\| \rho\|$ having a constant potential on $[a,b]$, and thus $k \mu + Bal(\rho,[a,b])=(k+\|\rho\|)\omega_{[a,b]}$. We apply the explicit formula for balayage onto an interval given in \cite[Eqn.\,(II.4.47)]{Saff-Totik}, and get for $t\in [a,b]$
\begin{align*}
   g(t) &:= \frac{d\mu}{d\omega_{[a,b]}}(t)
    = \frac{k+\| \rho \|}{k} - \frac{1}{k}
     \int_{-\infty}^a \frac{\sqrt{(b-y)(a-y)}}{t-y}\mathrm{d}\rho(y).
\end{align*}
As a consequence
$$
    0 \leq \lim_{t\to a+0} g(t) = \frac{k+\| \rho \|-\eta(a)}{k}  ,
$$
showing that $\eta(a)\leq k+\|\rho\|$ is finite. In particular, in case $a=\alpha$ we get \eqref{claim_alpha2}.
Moreover, by \cite[Eqn.\ (II.5.4)]{Saff-Totik}, with a suitable $F\in \mathbb R$,
\begin{align*}
    k(F-U^\mu(x)-Q(x)) &  = kF - U^{(k+\|\rho\|)\omega_{[a,b]}+\rho-Bal(\rho,[a,b])}(x)
    \\&= (k+\| \rho\|) g(x,\infty)
    - \int g(x,y) \, \mathrm{d}\rho(y),
\end{align*}
the right-hand side vanishing on $[a,b]$, and thus the constant $F$ coincides with the one in \eqref{eq.extremal2}. Hence, the above expression equals $kf(x)$, showing \eqref{eq.extremal5}.

It remains to show that Theorem~\ref{theorem-potential} implies
\eqref{eq.extremal6}, where we start to verify the hypotheses on
\begin{align*}
   g(t)
     &= \frac{d\mu}{d\omega_{[a,b]}}(t) =
     \frac{k+\| \rho \|-\eta(a)}{k}
     + \frac{t-a}{k} \int_{-\infty}^a
     \sqrt{\frac{b-y}{a-y}}\frac{\mathrm{d}\rho(y)}{t-y}.
\end{align*}
We first observe that $g$ is differentiable on $(a,b]$, with derivative
\[
      g'(t)= \frac{1}{k} \int_{-\infty}^a
     \sqrt{\frac{b-y}{a-y}}\frac{a-y}{(t-y)^2} \mathrm{d}\rho(y),
\]
which is both $\geq 0$ and decreasing in $t\in (a,b]$. Hence $g$ is increasing and concave in $(a,b)$, and, by a similar argument, $h(t):=g(t)/(t-a)$ is convex on $(a,b)$. Thus the assumptions of Theorem~\ref{theorem-potential} hold. With $P_k\in \Pi_k$ as in Theorem~\ref{theorem-potential} we have that
\begin{align*}
   \log \| w^k P_k \|_{[a,b]}
   &
   = \max_{x\in [a,b]} - k Q(x) + \log| P_k(x)|
   \\& \leq \max_{x\in [a,b]} - k Q(x) - kU^\mu(x) + C_{BW}
   = -k F + C_{BW} ,
\end{align*}
where for obtaining the inequality we have applied Theorem~\ref{theorem-potential}(a), and in the last equality we have used \eqref{eq.extremal2} and in particular the fact that $f$ vanishes on $[a,b]$. Also, for $x\in \mathbb R \setminus (a,b)$, we deduce from Theorem~\ref{theorem-potential}(b) and \eqref{eq.extremal2} that
\[
     \log w(x)^k |P_k(x)| \geq - k Q(x) - k U^\mu(x) = k f(x)-kF .
\]
Combining these two inequalities gives \eqref{eq.extremal6}.

\subsection{Proof of Theorem~\ref{Th_super}}\label{sec_proof_th_super}

For our proof of Theorem~\ref{Th_super}(a), we choose $n>d+1$ as in the statement such that \eqref{eq1.Th_super} and \eqref{eq2.Th_super} hold. By our assumption \eqref{eq1.Th_super} on $\lambda_j$, we may choose $\widetilde \lambda_j\leq \lambda_j$ such that
$$
     \sigma((-\infty,\widetilde \lambda_j])=\frac{j}{N}
     \quad \mbox{for $j=1,2,...,d+1$.}
$$
Consider $k=n-d\geq 2$, and the two measures of mass $d$
$$
    \rho = \delta_{\lambda_1}+...+\delta_{\lambda_d} , \quad
    \widetilde \rho = N \sigma|_{(-\infty,\widetilde \lambda_d]}.
$$
It is not too difficult to check that $U^{\widetilde \rho}(x)-U^\rho(x)$
is $\leq 0$ for $x\in [\lambda_{d+1},b]$, and $\geq 0$ for $x=0$. Hence, going back to the definition of $E_n([\lambda_{d+1},b])$, we get the chain of inequalities
\begin{eqnarray*}
   E_n([\lambda_{d+1},b]) &=& \min_{p\in \Pi_k} \frac{\| e^{-U^\rho} p \|_{[\lambda_{d+1},b]}}{e^{-U^\rho(0)} |p(0)|}
   \leq \min_{p\in \Pi_k} \frac{\| e^{-U^{\widetilde \rho}} p \|_{[\lambda_{d+1},b]}}{e^{-U^{\widetilde \rho}(0)} |p(0)|}
   \\&\leq& \exp\Bigl( C_{BW} + k U^\mu(0) + U^{\widetilde \rho}(0) - kF \Bigr),
\end{eqnarray*}
where in the last inequality we have applied Theorem~\ref{th_main} with $\alpha=\lambda_{d+1}$, and the external field $Q(x)=U^{\widetilde \rho/k}(x)$, and where the extremal measure $\mu$ 
and the constant $F$ are as in \eqref{eq.extremal1}, \eqref{eq.extremal2}.
On the other hand, with $t=n/N$, we know from \cite[Theorem~2.1]{Beck-CG} that
$$
      \exp\Bigl(C_{BW} - N \int_0^{n/N} g_{S(\tau)}(0,\infty) d\tau \Bigr)
      = \exp\Bigl(C_{BW} + n(U^{\nu_{t,\sigma}}(0) - C_{t,\sigma})\Bigr),
$$
with $\nu_{t,\sigma}$ the solution of the constrained equilibrium problem mentioned in the paragraph after \eqref{eq.BB_ABK}, and $C_{t,\sigma}$ the corresponding constant. Thus, for establishing Theorem~\ref{Th_super}(a), it only remains to show the inequality
\begin{equation} \label{super_claim}
      k U^\mu(x) + U^{\widetilde \rho}(x) - kF \leq
      n(U^{\nu_{t,\sigma}}(x) - C_{t,\sigma})
\end{equation}
for $x=0$.

Let us first show that \eqref{super_claim} holds for $x\in \supp(\mu)$.
Indeed, since $\supp(\mu)\subset [\lambda_{d+1},b]\subset [a(t),b]=\supp(\sigma/t-\nu_{t,\sigma})$ by assumption \eqref{eq2.Th_super},
we find from the respective equilibrium conditions for both extremal problems that both expressions on the left-hand side and on the right-hand side of
\eqref{super_claim} vanish for $x\in \supp(\mu)$. We also know that all measures involved in \eqref{super_claim} have finite energy,
with masses $\| k \mu +\widetilde \rho\| =k+d=n= n \, \| \nu_{t,\sigma} \|$. Let us show that $\widetilde \rho \leq n\nu_{t,\sigma}$.
Indeed, $\widetilde \lambda_d\leq \lambda_d<a(t)$ by construction and \eqref{eq2.Th_super}, and thus, by definition of $S(t)=[a(t),b]=\supp(\sigma/t-\nu_{t,\sigma})$,
$$
     n \nu_{t,\sigma}|_{(-\infty,\widetilde \lambda_d]}
     = \frac{n}{t} \sigma|_{(-\infty,\widetilde \lambda_d]}
     = N \sigma|_{(-\infty,\widetilde \lambda_d]} = \widetilde \rho .
$$
Hence, by subtracting $U^{\widetilde \rho}(x)$ from both sides of      \eqref{super_claim}, we get from the principle of domination for logarithmic potentials \cite[Theorem~II.3.2]{Saff-Totik} that \eqref{super_claim} holds for all $x\in \mathbb C$, and in particular for $x=0$, which concludes our proof of
Theorem~\ref{Th_super}(a).

For our proof of Theorem~\ref{Th_super}(b), we first observe that our assumption of Lemma~\ref{lemma_CG} on $\sigma$ and the choice of the eigenvalues $\lambda_{1,N}<\lambda_{2,N}<...$ of $A_N$, allows to show that not only \eqref{weak_limit} but also the quite technical \cite[Conditions (i)--(iv)]{Beck-CG} hold, we omit details. As a consequence of \cite[Theorem~2.2]{Beck-CG},
$$
    \lim_{{n,N\to \infty}\atop{n/N\to t}} E_n(\Lambda(A_N))^{1/n}
    =\exp\Bigl( - \frac{1}{t} \int_0^t g_{S(\tau)}(0,\infty)  d\tau\Bigr),
$$
that is, we have equality in \eqref{eq.BB_ABK}. Then, using Theorem~\ref{Th_super}(a) and the simple inequality $E_n(\Lambda(A_N) \leq E_n([\lambda_{d+1},b])$,
\begin{eqnarray*}
&&
   \liminf_{{n,N\to \infty}\atop {n/N\to t}} E_{n,N}([\lambda_{d_{n,N}+1,N},b])^{1/n}
   \leq \limsup_{{n,N\to \infty}\atop {n/N\to t}} E_{n,N}([\lambda_{d_{n,N}+1,N},b])^{1/n}
\\&&   \leq \exp\Bigl( - \frac{1}{t} \int_0^t g_{S(\tau)}(0,\infty)  d\tau\Bigr)
           =
   \lim_{{n,N\to \infty}\atop{n/N\to t}} E_n(\Lambda(A_N))^{1/n}
\\&&    \leq
           \liminf_{{n,N\to \infty}\atop {n/N\to t \atop 0\leq d <n-1}} E_{n,N}([\lambda_{d+1,N},b])^{1/n}.
    \leq
    \liminf_{{n,N\to \infty}\atop {n/N\to t}} E_{n,N}([\lambda_{d_{n,N}+1,N},b])^{1/n},
\end{eqnarray*}
and the statement of Theorem~\ref{Th_super}(b) follows.



\section{Conclusions}

In the particular case of an external field being given by the logarithmic potential of some positive measure supported on the left of $\Sigma$, we have shown that the weighted Bernstein-Walsh inequality is sharp up to some new universal constant $C_{BW}$. Our main tool is a variation of the technique of Totik of discretizing a logarithmic potential, provided that the underlying measure has a weight satisfying some monotonicity and/or convexity assumptions.

This new sharpness result for the weighted Bernstein-Walsh inequality leads to a variety of new explicit bounds for the convergence of conjugate gradients if we fix in advance a fixed number of small eigenvalues being considered as outliers. By approximately optimizing the number of outliers, we are able to partly show a conjecture formulated by Beckermann \& Kuijlaars \cite{Beck-CG} in terms of means of Green functions, and establish a new upper bound for conjugate gradients in form of an inequality for every iteration index $n$. Such bounds are of practical interest since the results of \cite{Beck-CG} are only of asymptotic nature. In addition, our bounds are valid for a single matrix and do no longer require to consider sequences of systems of equations with a joint eigenvalue distribution. We also give some (academic) numerical examples showing that this new bound perfectly matches the CG error (up to the choice of $C_{BW}$).

We believe that, with an optimal choice of $C_{BW}$, the quantity $e^{C_{BW}}$ is of modest size. This is clearly not true for our present explicit upper bound of $C_{BW}$, and remains a direction of future research, maybe asymptotic analysis could be helpful.

We also believe that our result on the discretization of a potential can be generalized to more general measures, for example without the assumption that
$t \mapsto \frac{g(t)}{t-a}$ is convex on $(a,b)$, which is used only once. This possibly would allow us to consider both small and large eigenvalues as outliers.

Finally, the above-mentioned conjecture on the CG convergence remains open for general sets $S(t)$.



\appendix


\section{Proof of the mean value property of Theorem~\ref{Theorem:mean:property}}\label{annexeA}


As said before, a central role in our analysis is played by the mean value property of the cumulative distribution function $W_g$ stated in Theorem~\ref{Theorem:mean:property}:
there exist constants $c_1 = \frac{1}{4}$ and $c_2 = \pi\sqrt{2}$ such that, for all $x,t\in [-1,1]$,
\begin{equation*}
   c_1 W_g'(\frac{t+x}{2}) \leq \frac{W_g(t)-W_g(x)}{t-x}  \leq c_2 W_g'(\frac{t+x}{2}) .
\end{equation*}
The aim of this section is to provide a proof of this mean value property.
We will first consider the two particular cases $g=1$ in Lemma~\ref{lemma:mean:property:Chebyshev:case} and $g(t)=1+t$ in Lemma~\ref{lemma:mean:property:droite case}.
The general case then will follow by concavity of $g$ and by convexity of $h(t)=g(t)/(t+1)$. In what follows it will be convenient to consider the substitution $t=-\cos(\alpha)$ and $x=-\cos(\beta)$, $\alpha,\beta \in [0,\pi]$, where we can suppose without loss of generality that $t> x$, and thus $0 \leq \beta < \alpha \leq \pi$.


\begin{Lemma}
  \label{lemma:mean:property:Chebyshev:case}
  For every $x,t \in [-1,1]$, we have for $c_3={\pi}/{\sqrt{2}}$
  \begin{equation}
  \label{lemma:mean:property:Chebyshev:case-equation}
    W_1'(\frac{t+x}{2}) \leq \frac{W_1(t)-W_1(x)}{t-x}  \leq c_3 W_1'(\frac{t+x}{2}) .
  \end{equation}
\end{Lemma}

\dem{
  Elementary trigonometric formulas give
\begin{align*}
  \frac{W_1(t)-W_1(x)}{t-x} & = \frac{k}{\pi} \frac{\alpha-\beta}{\cos(\beta)-\cos(\alpha)}
  = \frac{k}{\pi} \frac{\frac{\alpha-\beta}{2}}{\sin(\frac{\alpha-\beta}{2})} \frac{1}{\sin(\frac{\alpha+\beta}{2})} .
\end{align*}
Observing that $\frac{\alpha-\beta}{2} \in [0,\frac{\pi}{2}]$ and thus
\[
       \sin(\frac{\alpha-\beta}{2}) \leq \frac{\alpha-\beta}{2} \leq \frac{\pi}{2} \sin(\frac{\alpha-\beta}{2}),
\] we deduce that
\begin{equation*}
  \frac{k}{\pi} \frac{1}{\sin(\frac{\alpha+\beta}{2})} \leq \frac{W_1(t)-W_1(x)}{t-x} \leq \frac{k}{2} \frac{1}{\sin(\frac{\alpha+\beta}{2})} .
\end{equation*}
Since
\begin{equation*}
  W_1'(\frac{t+x}{2})
  = \frac{k}{\pi} \frac{1}{\sqrt{1-\cos^2(\frac{\alpha+\beta}{2})\cos^2(\frac{\alpha-\beta}{2})}},
\end{equation*}
the left-hand inequality in \eqref{lemma:mean:property:Chebyshev:case-equation} immediately follows.

If $\frac{\alpha+\beta}{2} \leq \frac{\pi}{2}$, then $0 \leq \frac{\alpha-\beta}{2} \leq \frac{\alpha+\beta}{2} \leq \frac{\pi}{2}$.
If $\frac{\alpha+\beta}{2} \geq \frac{\pi}{2}$, then $0 \leq \frac{\alpha-\beta}{2} \leq \pi - \frac{\alpha+\beta}{2} \leq \frac{\pi}{2}$.
In both cases we find that
\[ 1-\cos^2(\frac{\alpha+\beta}{2})\cos^2(\frac{\alpha-\beta}{2})  \leq 1-\cos^4(\frac{\alpha+\beta}{2}) \leq 2 \left( 1-\cos^2(\frac{\alpha+\beta}{2}) \right) ,  \]
which implies the right-hand side of \eqref{lemma:mean:property:Chebyshev:case-equation}.
}



We now turn to the special case $g(y)=1+y$ where we only require one inequality for $W_g=W_{1+y}$.
\begin{Lemma}
\label{lemma:mean:property:droite case}
  For every $x,t \in [-1,1]$, we have for $c_4={1}/{2}$
  \begin{equation*}
    c_4 W_{1+y}'(\frac{x+t}{2}) \leq \frac{W_{1+y}(t)-W_{1+y}(x)}{t-x}  .
  \end{equation*}
\end{Lemma}

\dem{
By Lemma~\ref{lemma:mean:property:Chebyshev:case},
\begin{align*}
 \frac{W_{1+y}(t)-W_{1+y}(x)}{t-x}
 &
 \geq
 \frac{W_{1+y}(t)-W_{1+y}(x)}{W_1(t)-W_1(x)}  W_1'(\frac{x+t}{2}) .
\end{align*}
 Thus it is sufficient to show that
\[ \frac{W_{1+y}(t)-W_{1+y}(x)}{W_1(t)-W_1(x)}  \geq \frac{1}{2} (1+\frac{x+t}{2}) . \]
 By definition of $W_{1+y}$,
\begin{align*}
  \frac{W_{1+y}(t)-W_{1+y}(x)}{W_1(t)-W_1(x)} & = \frac{1}{W_1(t)-W_1(x)} \int_{x}^{t} W_{1+y}'(s) \mathrm{d}s  \\
  & = \frac{ \alpha-\beta + \sin(\beta)-\sin(\alpha) }{\alpha-\beta} \\
  & = 1 - \frac{2}{\alpha-\beta} \sin(\frac{\alpha-\beta}{2}) \cos(\frac{\alpha+\beta}{2}) .
\end{align*}
Hence it remains to show that
\[
   \cos(\frac{\alpha+\beta}{2}) \left( 2 \frac{\sin(\frac{\alpha-\beta}{2})}{\frac{\alpha-\beta}{2}} - \cos(\frac{\alpha-\beta}{2}) \right) \leq 1 .
\]
Since $\gamma \mapsto 2 \sin(\gamma)-\gamma \cos(\gamma)$ is increasing in $[0,\pi/2]$, the factor in large brackets is $\geq 0$, and $\cos((\alpha+\beta)/2)\leq \cos((\alpha-\beta)/2$. Thus we only have to consider the worst case $\gamma=(\alpha+\beta)/2=(\alpha-\beta)/2\in [0,\pi/2]$, with
\[
   \cos (\gamma) \left( 2 \frac{\sin(\gamma)}{\gamma}- \cos(\gamma) \right)
   \leq 2\cos (\gamma) - \cos^2(\gamma)\leq 1.
\]
}


We are now prepared to give a proof of Theorem~\ref{Theorem:mean:property}.
For the upper bound, we use Lemma~\ref{lemma:mean:property:Chebyshev:case} in order to conclude that
\begin{align*}
  \frac{W_g(t)-W_g(x)}{t-x}
  & = \frac{W_g(t)-W_g(x)}{W_1(t)-W_1(x)} \frac{W_1(t)-W_1(x)}{t-x}
  \\
  & \leq \frac{W_g(t)-W_g(x)}{W_1(t)-W_1(x)} \, c_3 \, W_1'(\frac{t+x}{2}).
\end{align*}
Recalling that $g$ is concave, we get from the Jensen inequality
\begin{align*}
  \frac{W_g(t)-W_g(x)}{W_1(t)-W_1(x)} & = \int_{x}^{t} \frac{W_g'(s)}{W_1(t)-W_1(x)} \mathrm{d}s \\
  & = \int_{x}^{t} g(s) \frac{W_1'(s)}{W_1(t)-W_1(x)} \mathrm{d}s \\
  & \leq g \left( \int_{x}^{t} s \frac{W_1'(s)}{W_1(t)-W_1(x)} \mathrm{d}s \right)
  \\
  & \leq 2g(\frac{t+x}{2}) ,
\end{align*}
the last inequality being established in Lemma~\ref{lemma:u:moyenne:g} below.
Thus we obtain the upper bound with $c_2=2c_3= \pi \sqrt{2}$. For the lower bound, our argument is similar, but now we use Lemma~\ref{lemma:mean:property:droite case} in order to get
\begin{align*}
  \frac{W_g(t)-W_g(x)}{t-x}
  & = \frac{W_g(t)-W_g(x)}{W_{1+y}(t)-W_{1+y}(x)} \frac{W_{1+y}(t)-W_{1+y}(x)}{t-x}
  \\& \geq \frac{W_g(t)-W_g(x)}{W_{1+y}(t)-W_{1+y}(x)}
  c_4 W_{1+y}'(\frac{x+t}{2}).
\end{align*}
Recalling that $h(y)=g(y)/(1+y)$ is convex, we get from the Jensen inequality
\begin{align*}
  \frac{W_g(t)-W_g(x)}{W_{1+y}(t)-W_{1+y}(x)} & = \int_{x}^{t} h(s) \frac{W_{1+y}'(s)}{W_{1+y}(t)-W_{1+y}(x)} \mathrm{d}s \\
  & \geq h \left( \int_{x}^{t} s \frac{W_{1+y}'(s)}{W_{1+y}(t)-W_{1+y}(x)} \mathrm{d}s \right)
  \\
  & \geq \frac{1}{2} \frac{g(\frac{t+x}{2})}{1+\frac{t+x}{2}} ,
\end{align*}
where for the last inequality we apply Lemma~\ref{lemma:v:moyenne:g} below. This gives us the lower bound with $c_1=c_4/2= 1/4$.

For concluding, it remains to establish two technical results.



\begin{Lemma}
  \label{lemma:u:moyenne:g}
  For $-1\leq x < t \leq 1$ we have
  \[ g \left( \int_{x}^{t} \frac{s W_1'(s)}{W_1(t)-W_1(x)} \mathrm{d}s \right) \leq 2g(\frac{x+t}{2}). \]
\end{Lemma}

\dem{Elementary trigonometric computations give
\[ \overline x := \int_{x}^{t} \frac{s W_1'(s)}{W_1(t)-W_1(x)} \mathrm{d}s = -\cos \left( \frac{\alpha+\beta}{2}\right)  \frac{\sin \left( \frac{\alpha-\beta}{2} \right)}{\frac{\alpha-\beta}{2}} \in [-1,1],  \]
and
\[ \frac{x+t}{2} = -\cos \left( \frac{\alpha+\beta}{2} \right)  \cos \left( \frac{\alpha-\beta}{2} \right) . \]
We now have to distinguish two cases.
If $\frac{x+t}{2}\leq 0$ or, equivalently, $\cos( \frac{\alpha+\beta}{2})\geq 0$ then, using the fact that
$\sin \left( y \right) \geq y\cos \left( y \right) \geq 0$ for $y \in [0,\frac{\pi}{2}]$, we get $\overline x\leq \frac{x+t}{2}$,
and the statement $g(\overline x) \leq g(\frac{x+t}{2})$ follows (without a factor $2$) by monotonicity of $g$.
If however $\frac{x+t}{2}\geq 0$, then $ 1 \geq \overline x \geq \frac{x+t}{2} \geq 0$
by the same argument as in the first case.
By concavity of $g$ and \eqref{eq.tangent} we get
\begin{align*}
   g(\overline x) &\leq g(\frac{x+t}{2}) + g'(\frac{x+t}{2}-0)  \left( \overline x - \frac{x+t}{2} \right)
   \\& \leq g(\frac{x+t}{2}) + g'(\frac{x+t}{2}-0)
   \leq 2 g(\frac{x+t}{2}),
\end{align*}
the last inequality being shown in Lemma~\ref{properties:g:2}(b) below. Thus Lemma~\ref{lemma:u:moyenne:g} holds.
}

\begin{Lemma}
  \label{lemma:v:moyenne:g}
  For $-1\leq x < t \leq 1$ we have
  \[ h \left( \int_{x}^{t} s \frac{W_{1+y}'(s)}{W_{1+y}(t)-W_{1+y}(x)} \mathrm{d}s \right) \geq \frac{1}{2} h(\frac{x+t}{2}).  \]
\end{Lemma}

\dem{
Let us first show that
\begin{equation}
  \label{equation:v:larger:mean}
  \overline x:= \int_{x}^{t} s \frac{W_{1+y}'(s)}{W_{1+y}(t)-W_{1+y}(x)} \mathrm{d}s \geq \frac{x+t}{2} .
\end{equation}
We write shorter $w(s) = \frac{W_{1+y}'(s)}{W_{1+y}(t)-W_{1+y}(x)}$ being increasing in $s$. Hence
\begin{align*}
  \overline x - \frac{x+t}{2} & = \int_{x}^{t} (s-\frac{x+t}{2}) w(s) \mathrm{d}s \\
  & = \int_{x}^{t} (s-\frac{x+t}{2}) (w(s)-w(\frac{x+t}{2})) \mathrm{d}s \geq 0,
\end{align*}
as claimed in \eqref{equation:v:larger:mean}. Also, by definition, $\overline x\leq t$, and thus
\begin{align*}
   h(\overline x) &=\frac{g(\overline x)}{\overline x+1} \geq \frac{g(\frac{x+t}{2})}{\overline x+1}
   \geq \frac{g(\frac{x+t}{2})}{t+1}
   \geq \frac{1}{2} \frac{g(\frac{x+t}{2})}{1+\frac{x+t}{2}} = \frac{h(\frac{x+t}{2})}{2}.
\end{align*}
}




\section{Some further technical lemmata}\label{annexeB}

After having established the mean value property of $W_g$ in \S\ref{annexeA},
we gather in this section all the other technical properties of the abscissa $t_j=-\cos(\alpha_j)$ needed in \S 3.


In the sequel of this section we always suppose the conditions on $k,g,h$ of Theorem~\ref{theorem-potential} to be true,
that is, $k$ is some integer $\geq 2$, $g$ is non-negative, increasing and concave, and 
$h(t)=g(t)/(t+1)$ is convex.

The first result summarizes some properties of the function $g$.
\begin{Lemma}
  \label{properties:g:2}
The following properties hold:
  \begin{enumerate}
	\item[(a)] $h(t)=\frac{g(t)}{1+t}$ is decreasing on $(-1,1]$;
    \item[(b)] $g'(t-0)(1+t)  \leq g(t)$ for $t\in (-1,1)$, and $g'(t-0) \leq g(t)$ for $t\in [0,1)$;
	\item[(c)] $g(1) \leq 2$;
	\item[(d)] $g(1)\geq g(0)=h(0) \geq 1$.
  \end{enumerate}
\end{Lemma}

\dem{
  We first recall that, by concavity of $g$ on $[-1,1]$, we have for all
  $-1\leq x_1 < x_2 < x_3\leq 1$ that
  \begin{equation} \label{eq.concave}
       \frac{g(x_3)-g(x_2)}{x_3-x_2} \leq
       \frac{g(x_3)-g(x_1)}{x_3-x_1} \leq
       \frac{g(x_2)-g(x_1)}{x_2-x_1}
       .
  \end{equation}
  Since $g(-1)\geq 0$, we may therefore write
  $$
         h(t) = \frac{g(t)-g(-1)}{t-(-1)} + \frac{g(-1)}{t+1}
  $$
  as a sum of two decreasing functions, implying (a). Passing to the limit in \eqref{eq.concave}, we also have that the directional derivatives
  $g'(x_2-0)$ and $g'(x_2+0)$ exist for all $x_2\in (-1,1)$, with
  \[
      \frac{g(x_3)-g(x_2)}{x_3-x_2} \leq g'(x_2+0)\leq g'(x_2-0)\leq \frac{g(x_2)-g(x_1)}{x_2-x_1},
  \]
  and in particular
  \begin{equation} \label{eq.tangent}
      g(x)\leq g(x_2)+g'(x_2-0)(x-x_2) \quad \mbox{for all $x\in [-1,1]$.}
  \end{equation}
  Setting $x=-1$, $x_2=t$ in \eqref{eq.tangent} leads to (b) since $g(-1)\geq 0$.
  Furthermore, using the concavity of $g$ and setting $x_2=0$ in \eqref{eq.tangent}, we get for all $t\in [-1,1]$ that
  $$
        g(-1)\frac{1-t}{2} + g(1)\frac{1+t}{2} \leq g(t) \leq g(0) + t g'(0-).
  $$
  Taking into account \eqref{distribution},
  multiplying by $W_1'(t)$ and integrating from $-1$ to $1$ gives
  $$
      k \frac{g(1)}{2} \leq k \frac{g(1)+g(-1)}{2} \leq k \leq k g(0),
  $$
  implying parts (c) and (d).
}

The following elementary lemma will be helpful in what follows.
\begin{Lemma}
  \label{lemma:rapport:I_j:sin:angles}
  For $\gamma\geq 0$ and $0\leq \delta \leq \theta \leq \pi/2$ there holds
  \[
       \frac{\sin(\gamma)}{\sin(\delta)} \leq \frac{\theta}{\sin(\theta)} \frac{\gamma}{\delta} \leq \frac{\pi}{2}\frac{\gamma}{\delta}.
  \]
\end{Lemma}
\dem{Since $x\mapsto x/\sin(x)$ is increasing in $[0,\pi/2]$, we have that
  \[
     \frac{\sin(\gamma)}{\sin(\delta)}  \leq \frac{\gamma}{\sin(\delta)}
     = \frac{\delta}{\sin(\delta)} \frac{\gamma}{\delta}
     \leq
     \frac{\theta}{\sin(\theta)} \frac{\gamma}{\delta} \leq \frac{\pi}{2}\frac{\gamma}{\delta}.
  \]
}

The following result tells us that the angles $\alpha_j$ defined by $t_j=-\cos(\alpha_j)$ for $j=0,1,...,k$ have a quite regular behavior.

\begin{Lemma}
  \label{lemma:sequence:alpha:decreasing:rapport:I_j:avant}
  {\bf (a)} The sequence $(\alpha_{j+1}-\alpha_j)_{0\leq j <k-1}$ is decreasing. \\ {\bf (b)} For $j\in \{ 1,...,k-1\}$ there holds
  \[
        \frac{\alpha_{j+1}-\alpha_j}{\alpha_j}  \leq \frac{1}{j}  ,
  \]
  {\bf (c)} For $j\in \{ 0,...,k-1\}$ we have
  \[
     \frac{\alpha_{j+1}-\alpha_j}{\alpha_{j+1}}  \leq \frac{1}{j+1} .
  \]
  {\bf (d)} For $j \in \{ 0,1,...,k-2\}$
  \[
        \frac{\alpha_{k}-\alpha_j}{\alpha_{k}-\alpha_{j+1}}  \leq 4.
  \]
\end{Lemma}
\dem{
  Using that $g$ is increasing, we get for $j\in \{0,...,k-1\}$
\begin{align*}
  1 &= \int_{t_j}^{t_{j+1}} W_g'(t) \mathrm{d}t
  = \frac{k}{\pi} \int_{t_j}^{t_{j+1}} \frac{g(t)}{\sqrt{1-t^2}} \mathrm{d}t   \\&
  \left\{\begin{array}{l}
  \displaystyle  \geq \frac{k}{\pi} g(t_j) \int_{t_j}^{t_{j+1}} \frac{\mathrm{d}t}{\sqrt{1-t^2}} = \frac{k}{\pi} g(t_j) (\alpha_{j+1}-\alpha_j) ,
  \\
  \displaystyle \leq \frac{k}{\pi} g(t_{j+1}) \int_{t_j}^{t_{j+1}} \frac{\mathrm{d}t}{\sqrt{1-t^2}} = \frac{k}{\pi} g(t_{j+1}) (\alpha_{j+1}-\alpha_j) ,
  \end{array}\right.
\end{align*}
  implying that
\begin{equation}
  \label{equation1:lemma1}
  \frac{\pi}{kg(t_{j+1})} \leq \alpha_{j+1}-\alpha_j \leq \frac{\pi}{k g(t_j)} .
\end{equation}
  Thus (a) holds. For a proof of (b), we apply (a) to conclude that, for $j \in \{ 1,2,...,k-1 \}$,
\begin{align*}
  \frac{\alpha_{j+1}-\alpha_j}{\alpha_j}  =
    \frac{\alpha_{j+1}-\alpha_j}{\alpha_j-\alpha_0}  & = \frac{ \alpha_{j+1}-\alpha_j }{ \sum_{p=0}^{j-1} \alpha_{p+1}-\alpha_p } \\
  & \leq \frac{ \alpha_{j+1}-\alpha_j }{ j (\alpha_{j}-\alpha_{j-1}) } \leq \frac{1}{j} .
\end{align*}
  A proof of part (c) follows the same lines, we omit details.
  Let us finally show (d).
  In case $k=2$, we know from \eqref{equation1:lemma1} and Lemma~\ref{properties:g:2}(c) that $\alpha_2-\alpha_1 \geq \frac{\pi}{2g(1)} \geq \pi/4$, implying (d).
  In case $k\geq 3$ we can write
  \begin{align*}
     \frac{\alpha_{k}-\alpha_j}{\alpha_{k}-\alpha_{j+1}}
     & =
     1 + \frac{\alpha_{j+1}-\alpha_j}{\alpha_{k}-\alpha_{j+1}}
     \\& =  1 + \frac{\alpha_{j+1}-\alpha_j}{\sum_{\ell=j+1}^{k-1}
     (\alpha_{\ell+1}-\alpha_{\ell})}
     \\&\leq 1 + \frac{1}{k-1-j} \frac{\alpha_{j+1}-\alpha_j}{\alpha_{k}-\alpha_{k-1}},
  \end{align*}
  where in the last inequality we have applied (a).
	By part (c), $\alpha_{j+1}-\alpha_j \leq \alpha_{j+1}/(j+1)\leq \pi/(j+1)$,
	and $\alpha_{k}-\alpha_{k-1}\geq \pi/(kg(t_k))\geq \pi/(2k)$ by \eqref{equation1:lemma1} and Lemma~\ref{properties:g:2}(c).
	Hence using that $k\geq 3$, we obtain
  $$
     \frac{\alpha_{k}-\alpha_j}{\alpha_{k}-\alpha_{j+1}}
     \leq 1 + 2 \frac{k}{(j+1)(k-1-j)} \leq 1 + 2 \frac{k}{k-1}\leq 4.
  $$
}

%

The following result is used in our proof of Proposition~\ref{general:case:sum:3}.
\begin{Lemma}
  \label{lemma:rapport:I_j:avant}
  For $j\in \{ 1,2,...,k-1\}$ there holds
  \[ 1+\frac{t_j+t_{j+1}}{2} \leq c_5 (1+t_j) \]
  where $c_5 = \frac{3\pi}{4}+1$.
\end{Lemma}
\dem{
  By Lemma~\ref{lemma:sequence:alpha:decreasing:rapport:I_j:avant}(b),
  $\alpha_{j+1}-\alpha_j \leq \alpha_j/j\leq \pi$, implying that
  \[
     \sin(\frac{\alpha_{j+1}-\alpha_j}{2}) \leq \sin(\frac{\alpha_j}{2j}) \leq \sin(\frac{\alpha_j}{2}) .
  \]
  Moreover, since $\alpha_{j+1}+\alpha_j \leq (2+1/j) \alpha_j$, we get by Lemma~\ref{lemma:rapport:I_j:sin:angles}
  \begin{align*}
     \frac{t_{j+1}-t_j}{t_{j}-t_0} & = \frac{ \sin(\frac{\alpha_{j+1}-\alpha_j}{2})}{\sin(\frac{\alpha_j}{2})}
     \frac{\sin(\frac{\alpha_{j+1}+\alpha_j}{2}) }{\sin(\frac{\alpha_j}{2})} \\
     & \leq \frac{\pi}{2} (2+\frac{1}{j}) \leq \frac{3 \pi}{2},
  \end{align*}
  and Lemma~\ref{lemma:rapport:I_j:avant} follows.
}

Let us now show the two main properties required for our proof of Proposition~\ref{general:case:sum:1}.

\begin{Lemma}
  \label{lemma:rapport:I_j:apres}
	For $j\in \{ 0,1,...,k-2\}$ we have
  \[ 1-\frac{t_j+t_{j+1}}{2} \leq 9 (1-t_{j+1}). \]
\end{Lemma}
\dem{
  We will show the equivalent statement
\begin{align*}
  t_{j+1}-t_j \leq 16 (t_{k}-t_{j+1}) .
\end{align*}
If $t_{j+1} \leq {1}/{\sqrt{2}}$, we obtain
$$
   \frac{t_{j+1}-t_j}{1-t_{j+1}} \leq \frac{t_{j+1}+1}{1-t_{j+1}}
   \leq (1+\sqrt{2})^2 \leq 6 \leq 16.
$$
 It remains to consider the case $t_{j+1} \geq {1}/{\sqrt{2}}$, and thus $\alpha_{j+1} \geq {3\pi}/{4}$, or $(\pi-\alpha_{j+1})/{2} \leq {\pi}/{8}$.
Using first Lemma~\ref{lemma:rapport:I_j:sin:angles} and then Lemma~\ref{lemma:sequence:alpha:decreasing:rapport:I_j:avant}(d),
we obtain
\begin{align*}
  \frac{t_{j+1}-t_j}{1-t_{j+1}} & = -1 + \frac{1-t_j}{1-t_{j+1}} = -1 + \frac{ 1+\cos(\alpha_j) }{1+\cos(\alpha_{j+1})} \\
	& = -1 + \frac{ \cos^2(\frac{\alpha_j}{2}) }{\cos^2(\frac{\alpha_{j+1}}{2})}
	= -1 + \frac{ \sin^2(\frac{\pi}{2}-\frac{\alpha_j}{2}) }{\sin^2(\frac{\pi}{2}-\frac{\alpha_{j+1}}{2})} \\
	& \leq -1 + \left( \frac{\pi/8}{\sin(\pi/8)}  \frac{ \pi - \alpha_j }{ \pi - \alpha_{j+1} } \right)^2
	 \leq -1 + 16 \left( \frac{\pi/8}{\sin(\pi/8)}\right)^2 \leq 16.
\end{align*}
}

\begin{Lemma}
  \label{lemma:rapport:I_j:dur}
  For $j\leq j_0-1\leq k-2$ we have
  \[ \frac{1+\frac{t_{j+1}+t_{j_0}}{2}}{1+t_{j+1}} \leq  \frac{\pi^2}{8} \frac{j_0^2}{(j+1)^2}.
  \]
\end{Lemma}
\dem{
  Notice that, by Lemma~\ref{lemma:rapport:I_j:sin:angles},
\begin{align*}
   \frac{1+\frac{t_{j+1}+t_{j_0}}{2}}{1+t_{j+1}} & = 1 + \frac{t_{j_0}-t_{j+1}}{2(1+t_{j+1})} =
	 1 + \frac{ \sin(\frac{\alpha_{j+1}+\alpha_{j_0}}{2}) \sin(\frac{\alpha_{j_0}-\alpha_{j+1}}{2}) }{2\sin^2(\frac{\alpha_{j+1}}{2})} \\
   & \leq 1 + \frac{\pi^2}{8} \frac{ (\alpha_{j+1}+\alpha_{j_0}) (\alpha_{j_0}-\alpha_{j+1}) }{\alpha_{j+1}^2}.
\end{align*}
   Applying Lemma~\ref{lemma:sequence:alpha:decreasing:rapport:I_j:avant}(a),
   and recalling that $j_0\geq j+1$, we obtain
\begin{align*}
  \frac{\alpha_{j_0}-\alpha_{j+1}}{\alpha_{j+1}}
  & = \frac{ \sum_{\ell=j+1}^{j_0-1} (\alpha_{\ell+1}-\alpha_\ell) }{ \sum_{\ell=0}^{j} (\alpha_{\ell+1}-\alpha_\ell) } \\
  & \leq \frac{(j_0-j-1) (\alpha_{j+2}-\alpha_{j+1}) }{ (j+1) (\alpha_{j+1}-\alpha_{j}) } \leq \frac{j_0-j-1}{j+1},
\end{align*}
and
\begin{align*}
  \frac{\alpha_{j_0}+\alpha_{j+1}}{\alpha_{j+1}}
  & = \frac{\alpha_{j_0}-\alpha_{j+1}}{\alpha_{j+1}} + 2
   \leq \frac{j_0+j+1}{j+1} .
\end{align*}
Combining the three inequalities, we deduce that
\begin{align*}
   \frac{1+\frac{t_{j+1}+t_{j_0}}{2}}{1+t_{j+1}}
   & \leq
   1 + \frac{\pi^2}{8} \frac{j_0^2-(j+1)^2}{(j+1)^2}
   \leq  \frac{\pi^2}{8} \frac{j_0^2}{(j+1)^2} .
\end{align*}
}

The three following results are required in our proof of Proposition~\ref{general:case:sum:2}.


\begin{Lemma}
  \label{lemma:general:case:position:t1}
  For $k\geq 2$, we have
  \[ t_{k-1} \geq 0 . \]
\end{Lemma}
\dem{
  Suppose that $t_{k-1}<0$.
	Then using Lemma~\ref{properties:g:2}(d),
	and the fact that $g$ is increasing allows us to find a contradiction
  \begin{align*}
    1 & = \int_{t_{k-1}}^{1} W_g'(t) \mathrm{d}t > \int_{0}^{1} W_g'(t) \mathrm{d}t = \int_{0}^{1} \frac{g(t) k}{\pi} \frac{1}{\sqrt{1-t^2}} \mathrm{d}t \\
    & \geq \frac{g(0) k}{\pi} \int_{0}^{1} \frac{1}{\sqrt{1-t^2}} \mathrm{d}t \geq \frac{k}{2}.
  \end{align*}
}

\begin{Lemma}
  \label{lemma:general:case:racine:t1-t0}
  There holds
	\begin{enumerate}
	  \item[(a)] $\sqrt{t_1-t_0} \leq \frac{3\pi}{\sqrt{2}g(t_1)k}$;
      \item[(b)] $\sqrt{t_k-t_{k-1}} \leq \frac{\pi}{kg(1)}$;
      \item[(c)] For all $j\in \{0,1,...,k-1\}$ we have $t_{j+1}-t_j \leq ( \frac{12\pi}{k} )^{1/3}$.
  \end{enumerate}
\end{Lemma}
\dem{
   By Lemma~\ref{properties:g:2}(a), we find for $t\in (t_0,t_1]$ that
   $$
      W_g'(t) = h(t)\frac{k}{\pi} \frac{\sqrt{1+t}}{\sqrt{1-t}} \geq \frac{h(t_1)k}{\pi\sqrt{2}} \sqrt{1+t}.
   $$
   Integrating over the interval $[t_0,t_1]=[-1,t_1]$ gives
   \[
      1 \geq \frac{h(t_1)k\sqrt{2}}{3\pi}(1+t_1)^{3/2} = \frac{g(t_1)k\sqrt{2}}{3\pi}(1+t_1)^{1/2} ,
   \]
   which implies part (a).
   By Lemma~\ref{properties:g:2}(a) and Lemma~\ref{lemma:general:case:position:t1}, there holds for $t\in [t_{k-1},t_k]\subset [0,1]$,
   \[
        W_g'(t) \geq \frac{k}{\pi} \frac{h(1)}{\sqrt{1-t}},
   \]
   and by integrating over the interval $[t_{k-1},t_k]=[t_{k-1},1]$
   we get
   \[
       1 \geq \frac{2kh(1)}{\pi} (1-t_{k-1})^{1/2} = \frac{kg(1)}{\pi} (1-t_{k-1})^{1/2} ,
   \]
   as required for part (b). For a proof of (c), we observe that,
   by Lemma~\ref{lemma:sequence:alpha:decreasing:rapport:I_j:avant}(a),
   \begin{align*}
       t_{j+1} - t_{j}
       & = 2 \sin(\frac{\alpha_{j+1}-\alpha_{j}}{2})
       \sin(\frac{\alpha_{j+1}+\alpha_{j}}{2})
       \\& \leq 2 \sin(\frac{\alpha_{1}-\alpha_{0}}{2})
       = \sqrt{2(t_{1}-t_{0})}.
   \end{align*}
   By concavity and positivity of $g$ and Lemma~\ref{properties:g:2}(d),
   \begin{align*}
      g(t_{1}) &
      \geq g(t_{0})\frac{1-t_{1}}{2}+ g(1) \frac{t_{1}-t_{0}}{2}
      \geq \frac{t_{1}-t_{0}}{2} .
   \end{align*}
   Multiplying with $\sqrt{t_{1}-t_{0}}$ and applying part (a) we arrive at
   \[
       (t_{1}-t_{0})^{3/2} \leq 2g(t_1) \sqrt{t_1-t_0} \leq \frac{3\pi\sqrt{2}}{k},
   \]
   which yields part (c).
}

\begin{Lemma}
  \label{lemma:rapport:I_j:I_j+1}
	For every $j\in\{0,\dots,k-2\}$
  \[  \frac{2}{3\pi} \leq \frac{t_{j+1}-t_{j}}{t_{j+2}-t_{j+1}} \leq \frac{6c_2\sqrt{c_5}}{c_1} . \]
\end{Lemma}
\dem{
  In order to show the left-hand inequality, we write
\begin{align*}
  \frac{t_{j+2}-t_{j+1}}{t_{j+1}-t_{j}} & = \frac{\sin(\frac{ \alpha_{j+1}+\alpha_{j+2} }{2} )}{\sin(\frac{ \alpha_j + \alpha_{j+1} }{2})}
  \frac{\sin(\frac{ \alpha_{j+2}-\alpha_{j+1} }{2} )}{\sin(\frac{ \alpha_{j+1}-\alpha_j }{2} )}
  \leq \frac{\sin(\frac{ \alpha_{j+1}+\alpha_{j+2} }{2} )}{\sin(\frac{ \alpha_j + \alpha_{j+1} }{2})},
\end{align*}
where we have applied Lemma~\ref{lemma:sequence:alpha:decreasing:rapport:I_j:avant}(a).
We claim that the right-hand term is $\leq 3\pi/2$.
Indeed, if ${ (\alpha_j + \alpha_{j+1}) }/{2} \geq \pi/2$, then this quotient is less than one. Else,
by using Lemma~\ref{lemma:rapport:I_j:sin:angles}, we obtain
$$
    \frac{\sin(\frac{ \alpha_{j+1}+\alpha_{j+2} }{2} )}{\sin(\frac{ \alpha_j + \alpha_{j+1} }{2} ) } \leq \frac{\pi}{2} \frac{\alpha_{j+1}+\alpha_{j+2}}{\alpha_j + \alpha_{j+1}  },
$$
and from Lemma~\ref{lemma:sequence:alpha:decreasing:rapport:I_j:avant}(b) we get that
$\alpha_{j+1}+\alpha_{j+2} = \alpha_{j+2}-\alpha_{j+1}+2\alpha_{j+1} \leq 3\alpha_{j+1} \leq 3 (\alpha_j + \alpha_{j+1})$.

To prove the right-hand inequality in Lemma~\ref{lemma:rapport:I_j:I_j+1}
we use Theorem~\ref{Theorem:mean:property} and Lemma~\ref{properties:g:2}(a) in order to obtain
\begin{align*}
  \frac{t_{j+1}-t_{j}}{t_{j+2}-t_{j+1}} & \leq \frac{c_2}{c_1} \frac{W_g'(\frac{t_{j+1}+t_{j+2}}{2})}{W_g'(\frac{t_{j}+t_{j+1}}{2})} \\
  & \leq \frac{c_2}{c_1} \frac{h(\frac{t_{j+1}+t_{j+2}}{2}) (1+\frac{t_{j+1}+t_{j+2}}{2}) }{h(\frac{t_{j}+t_{j+1}}{2})
	(1+\frac{t_{j}+t_{j+1}}{2})} \frac{W'(\frac{t_{j+1}+t_{j+2}}{2})}{W'(\frac{t_{j}+t_{j+1}}{2})} \\
  & \leq \frac{c_2}{c_1} \sqrt{\frac{1+\frac{t_{j+1}+t_{j+2}}{2} }{1+\frac{t_{j}+t_{j+1}}{2}}
	 \frac{1-\frac{t_{j}+t_{j+1}}{2}}{1-\frac{t_{j+1}+t_{j+2}}{2}}}.
\end{align*}
With help of Lemma~\ref{lemma:rapport:I_j:avant} and Lemma~\ref{lemma:rapport:I_j:apres} we obtain
\begin{align*}
  \frac{t_{j+1}-t_{j}}{t_{j+2}-t_{j+1}}
	& \leq \frac{c_2}{c_1} 3\sqrt{c_5} \sqrt{\underbrace{\frac{1+t_{j+1} }{1+\frac{t_{j}+t_{j+1}}{2}}}_{\leq 2}
	 \underbrace{\frac{1-t_{j+1}}{1-\frac{t_{j+1}+t_{j+2}}{2}}}_{\leq 2}}
    \leq \frac{6c_2\sqrt{c_5}}{c_1}
\end{align*}
}



\end{document}